\documentclass[12pt]{article}
\usepackage{epsfig,latexsym,amssymb}
\input amssym.def
\input amssym.tex
\newcommand{\eq}{\begin{equation}}
\newcommand{\en}{\end{equation}}
\newcommand{\tabo}{1}
\newcommand{\tabt}{2}

\newcommand{\Pin}{\Pi_{n,N}}

\newcommand{\gamhn}{\Gamma_{h,n}}
\newcommand{\bmbeta}{\beta}
\newcommand{\HH}{H}
\newcommand{\tX}{X^*}

\newcommand{\bmB}{B}

\newcommand{\re}[1]{\mbox{(\ref{#1})}}
\newcommand{\rem}[1]{\mbox{{\em (\ref{#1})}}}

\def\qed{\mbox{\rule{0.5em}{0.5em}}}

\newtheorem{Theorem}{Theorem}
\newtheorem{theorem}[Theorem]{Theorem}
\newtheorem{lemma}[Theorem]{Lemma}
\newtheorem{corollary}[Theorem]{Corollary}
\newtheorem{construction}[Theorem]{Construction}
\newtheorem{proposition}[Theorem]{Proposition}
\newtheorem{example}[Theorem]{Example}
\newtheorem{defn}[Theorem]{Definition}

\newtheorem{question}[Theorem]{Question}
\newtheorem{conjecture}[Theorem]{Conjecture}
\newtheorem{condition}[Theorem]{Condition}
\newtheorem{remark}[Theorem]{Remark}
\newtheorem{problem}[Theorem]{Problem}

\def\proof{\noindent{\bf Proof.\ \ }}
\def\note{\noindent{\bf Note.\ \ }}

\newcommand {\ints}  {\mathbb{Z}}
\newcommand {\reals} {\mathbb{R}}
\newcommand {\real} {\mathbb{R}}
\newcommand {\complex} {\mathbb{C}}

\newcommand {\sigo}{\Sigma_1}
\newcommand {\sigtd}{\Sigma_{2,d}}
\newcommand {\sigh}{\Sigma_h}
\newcommand {\sighp}{\Sigma_h^{\#}}
\newcommand {\sigoh}{\hat{\Sigma}_1}
\newcommand {\sigt}{\Sigma_2}
\newcommand {\ST}{\sigt}

\newcommand {\sigop}{\Sigma_1^{\#}}

\newcommand {\sigtp}{\Sigma_2^{\#}}

\newcommand{\te}{\rightarrow}
\newcommand{\ed}{\mbox{$ \ \stackrel{d}{=}$ }}
\newcommand{\convd}{\mbox{$ \ \stackrel{\!d}{\rightarrow}$ }}
\newcommand{\eps}{\varepsilon}

\newcommand{\expe}{\varepsilon}

\newcommand{\Lev}{L\'evy}
\newcommand{\LK}{L\'evy-Khintchine}

\newcommand{\Tak}{T{\'a}kacs}

\newcommand{\Pol}{P\'olya}
\newcommand{\Polya}{P\'olya}

\newcommand{\lb}[1]{\label{#1}}

\newenvironment{thm}[1]{\begin{theorem}\label{#1}}{\end{theorem}}

\newenvironment{lmm}[1]{\begin{lemma}\label{#1}}{\end{lemma}}

\newenvironment{prp}[1]{\begin{proposition}\protect\label{#1}}{\end{proposition}}

\begin{document}
\title{Probability laws related to the Jacobi theta and 
Riemann zeta functions, and Brownian excursions}
\author{Philippe Biane\thanks{CNRS, DMA, 45 rue d'Ulm
75005 Paris, France }, 
~Jim Pitman\thanks {Dept. Statistics, U. C., Berkeley.  
Research supported in part by N.S.F. Grant DMS97-03961. 
This paper is also available as Technical Report No. 569, U. C.
Berkeley}
~and ~Marc Yor\thanks{Laboratoire de Probabilit\'es, 
Universit\'e Pierre et Marie Curie, 4 Place Jussieu F-75252 
Paris Cedex 05, France.}}

\date{\today}
\maketitle
\vskip 1cm
\newcommand{\Tinf}[2]{ T_{\inf}^{#1,#2} }
\newcommand{\EP}{Euler Process}
\newcommand{\XX}{X} 
\newcommand{\vph}{\varphi} 
\newcommand{\tM}{\tilde{M}} 
\newcommand{\hF}{\tilde{F}} 
\newcommand{\tV}{\tilde{V}} 
\newcommand{\sqsig} {\mbox{$\sqrt{ {\pi \over 2 } \ST }$}}
\newcommand{\sqpt} {\mbox{$\sqrt{ {\pi \over 2 } t}\,$}}
\newcommand{\sqpv} {\mbox{$\sqrt{ {\pi \over 2 } v}$}}
\newcommand{\fsig}{f_{\ST}} 
\newcommand{\lt}{\varphi} 
\newcommand{\Tex}{R}
\newcommand{\epx}{\varepsilon}
\newcommand{\nuf}{\bar{\nu}}
\newcommand{\Tsup}[2]{ T_{\sup}^{#1,#2} }
\newcommand {\besd}{ \mbox{BES$(d)$}}
\newcommand {\bes}{ \mbox{BES}}
\newcommand {\etil}{ \hat{\eta}_1 }
\newcommand {\besdbr}{ \mbox{BES$(\delta)^{\br}$}}
\newcommand {\besfd}{ \mbox{BES$(4-\delta)$}}
\newcommand {\best}{ \mbox{BES$(3)$}}
\newcommand {\BEST}{\mbox{BES$_0(3)$}}
\newcommand {\TT}{\Sigma}
\newcommand {\CT}{\tilde{\Sigma}}
\newcommand {\bestbr}{ \mbox{BES$(3)^{\br}$}}
\newcommand {\Pbly}{ P^\dagger }
\newcommand {\Pbr}{ P^{\rm br}}
\newcommand {\Ed}{ P^{(\delta)}}
\newcommand {\ee}{e}
\newcommand {\giv}{ \,|\,}
\newcommand {\epo}{ \varepsilon_0}
\newcommand {\epso}{ \varepsilon}
\newcommand {\epone}{ \varepsilon_1}
\newcommand {\eptwo}{ \varepsilon_2}
\newcommand {\epj}{ \varepsilon_j}
\newcommand {\epi}{ \varepsilon_i}
\newcommand {\epsi}{ \tilde{\varepsilon}_k}
\newcommand {\tjs}{ T_{j}^*}
\newcommand {\tjsm}{ T_{j-1}^*}
\newcommand {\tks}{ T_{k}^*}
\newcommand {\tksm}{ T_{k-1}^*}
\newcommand{\bX}{{\bf X}}
\newcommand{\kl}{{\lambda}}
\newcommand{\Te}{{T_\kl}}
\newcommand{\GT}{{G_{T_\la}}}
\newcommand{\FF}{{\cal F}}
\newcommand{\var}{{\rm Var}}
\font\bb=msbm10
\def\bR{\hbox{\bb R}}
\newcommand {\BLY}{ \mbox{$B^\dagger$}}
\newcommand{\br}{\mbox{$\scriptstyle{\rm br}$}}
\newcommand{\ex}{\mbox{$\scriptstyle{\rm ex}$}}
\newcommand{\me}{\mbox{$\scriptstyle{\rm me}$}}
\newcommand{\birth}{\mbox{$\scriptstyle{\rm birth}$}}
\newcommand{\BB}{B^{\br}}
\newcommand{\RB}{R}
\newcommand{\Bex}{B^{\ex}}
\newcommand{\Mex}{M^{\ex}}
\newcommand{\Bme}{B^{\me}}
\newcommand{\exj}{B^{\ex}_j}
\newcommand{\tb}{T^{(3)}}
\newcommand{\ts}{T^{*}}
\newcommand{\la}{\lambda}
\newcommand{\emu}{\lambda}
\newcommand{\csch}{{\rm csch} }
\newcommand{\sech}{{\rm sech} }
\newcommand{\ins}{\nu}
\newcommand{\nuF}{\nu_F^{\br}}
\newcommand{\nuV}{\nu_V^{\br}}
\newcommand{\nb}{N^{\birth}}
\newcommand{\nst}{N^{*}}
\newcommand{\nstm}{N^{*}(t-)}
\newcommand{\muinv}{\mu^{-1}}
\newcommand{\muu}{\mu}
\newcommand{\ggla}{\la^{-\gamma} \Gamma_{\alpha}}
\newcommand{\gamal}{\Gamma_{\alpha}}
\newcommand{\gama}{\Gamma_{a}}
\newcommand{\gamh}{\Gamma_{h}}
\newcommand{\hal}{H_\alpha}
\newcommand{\hnn}{H_\nn}
\newcommand{\zetam}{\zeta_{-}}
\newcommand{\tho}{F_1}
\newcommand{\sumnn}{\sum_{n = - \infty}^\infty}
\newcommand{\bab}{B_{a,b}}
\newcommand{\hnnp}{H_{\nu}'}
\newcommand{\hnno}{\hat{H}_\nu}
\newcommand{\nn}{\nu}
\newcommand{\halp}{H_{\alpha}'}
\newcommand{\al}{\alpha}
\newcommand{\ca}{C_\alpha}
\newcommand{\hhalf}{H_{1/2}}
\newcommand{\hmhalf}{H_{-1/2}}
\newcommand{\kal}{k_\alpha}
\newcommand{\calph}{c_\alpha}
\newcommand{\Pal}{P_{2 - 2 \alpha}}
\newcommand{\Palt}{\widetilde{P}_\alpha}
\newcommand{\Eal}{E_\alpha}
\newcommand{\gamj}{\Gamma_j}
\newcommand{\overl}[1]{\overline{#1}}
\newcommand{\underl}[1]{\underline{#1}}
\newcommand{\hf}{ \mbox{${1 \over 2}$}}
\newcommand{\thf}{ \mbox{${3 \over 2}$}}
\newcommand{\qr}{ \mbox{${1 \over 4}$}}
\newcommand{\tqr}{ \mbox{${3 \over 4}$}}
\newcommand{\hfs}{ \mbox{${s \over 2}$}}
\newcommand{\rtop}{ \mbox{$\sqrt{2 \over \pi}$}}
\renewcommand{\top}{ \mbox{${2 \over \pi}$}}
\newcommand{\rpot}{ \mbox{$\sqrt{\pi \over 2}$}}
\newcommand{\pot}{\mbox{${\pi \over 2}$}}
\newcommand{\hl}{ - \mbox{${1 \over 2}$} \lambda^2}
\newcommand{\hfpi}{ \mbox{${\pi \over 2}$}  }
\newcommand{\BBt}{B_t^{\br}}
%%%newdefs
\newcommand{\mbrd}{M_\delta} 
\newcommand{\mbr}{M} 
\newcommand{\mex}{M^{\ex}}
\newcommand{\mbrone}{M_1} 
\newcommand{\Pdel}{P_\delta} 
\newcommand{\mexd}{M_\delta^{\ex}} 
\newcommand{\mbral}{M_{2 - 2 \alpha}}
\newcommand{\hmal}{H_{-\alpha}} 
\newcommand{\mbrpj}{M_j^{\br +}}
\newcommand{\mbrpk}{M_k^{\br +}}
\newcommand{\mbrpo}{M_1^{\br +}}
\newcommand{\mbrnj}{M_j^{\br -}}
\newcommand{\mbrno}{M_1^{\br -}}
\newcommand{\mbrj}{M_j^{\br}}
\newcommand{\fbrj}{F_j^{\br}}
\newcommand{\tfbrj}{\widetilde{F}_j^{\br}}
\newcommand{\fbrjh}{\widehat{F}_j^{\br}}
\newcommand{\fbrjhm}{\widehat{F}_{j-1}^{\br}}
\newcommand{\fj}{F_j}
\newcommand{\vbrj}{V_j^{\br}}
\newcommand{\fbrjo}{F^{\br}_1}
\newcommand{\fbrjt}{F^{\br}_2}
\newcommand{\vbro}{V^{\br}_1}
\newcommand{\vbrt}{V^{\br}_2}
\newcommand{\vbrth}{V^{\br}_3}
\newcommand{\mbrp}{M^{{\br} +}}
\newcommand{\mbrn}{M^{{\br} -}}
\newcommand{\mbro}{M_1^{\br}}
\newcommand{\mbrt}{M_2^{\br}}
\newcommand{\fbro}{F_1^{\br}}
\newcommand{\fbrt}{F_2^{\br}}
\newcommand{\mbrpt}{M_2^{\br + }}
\newcommand{\mbrnt}{M_2^{\br - }}
\newcommand{\mbrnseq}{(M_j^{\br - })}
\newcommand{\mbrpseq}{(M_j^{\br + })}
\newcommand{\mbrseq}{(M_j^{\br})}
\newcommand{\tfbrseq}{(\widetilde{F}_j^{\br})}
\newcommand{\vbrseq}{(V_j^{\br})}
\newcommand{\fbrseq}{(F_j^{\br})}
\newcommand{\mbrseqst}{(M_j^{*})}
\newcommand{\mbrjst}{M_j^{*}}
\def\qed{\mbox{\rule{0.5em}{0.5em}}}
\def\proof{\noindent{\bf Proof.\ \ }}
\def\note{\noindent{\bf Note.\ \ }}
\def\endpf{\hfill $\Box$ \vskip .25in}

\begin{abstract}
This paper reviews known results which 
connect Riemann's integral representations of his zeta function,
involving Jacobi's theta function and its derivatives,
to some particular probability laws governing sums
of independent exponential variables. These laws are related
to one-dimensional Brownian motion and  
to higher dimensional Bessel processes. We present some 
characterizations of these probability laws, and some
approximations of Riemann's zeta function which are related to these laws.
\end{abstract}


\newpage
\newpage
\tableofcontents
\newpage

\section{Introduction}

In his fundamental paper \cite{riemann59}, Riemann showed that
the Riemann zeta function, initially defined by the series 
\eq
\lb{zeta1}
\zeta(s) :=\sum_{n=1}^{\infty}n^{-s}
~~~~~~~~ (\Re s>1)
\en
admits a meromorphic
continuation to the entire complex plane, with only a simple pole 
at 1,
and that the function
\eq
\lb{xifn}
\xi(s):= \hf s ( s-1) \pi^{-s/2} \Gamma( \hf s ) \zeta( s)  ~~~~~~~~(\Re s  > 1)
\en
is the restriction to $(\Re s  > 1)$ of a unique entire analytic function
$\xi$, which satisfies the functional equation
\eq
\lb{xieq}
\xi(s)=  \xi (1-s) 
\en
for all complex $s$. 
These basic properties of $\zeta$ and 
$\xi$ follow from a representation of $ 2 \xi$ as the Mellin transform of a
function involving derivatives of Jacobi's theta function. 
This function turns out to be the density of a probability distribution on 
the real line, which has deep and intriguing connections with the theory of Brownian motion.
This distribution first appears in the probabilistic 
literature in the 1950's in the work of Feller \cite{fel51}, 
Gnedenko \cite{gned54}, and \Tak\  \cite{tak58},
who derived it as the asymptotic distribution as $n \te \infty$
of the range of a simple one-dimensional random walk conditioned to 
return to its origin after $2n$ steps, and found formula \re{y1} below
for $s = 1,2, \cdots$.
Combined with the approximation of random walks by Brownian motion,
justified by Donsker's theorem \cite{bil68,durrett95,ry99},
the random walk asymptotics imply that if
\eq
\lb{bmm}
Y := \rtop \left( \max_{0 \le u \le 1} b_u - \min_{0 \le u \le 1} b_u \right)
\en
where
$(b_u, 0 \le u \le 1)$ is the standard Brownian bridge derived 
by conditioning a one-dimensional Brownian motion $(B_u, 0 \le u \le 1)$
on $B_0 = B_1 = 0$, then 
\eq
\lb{y1}
E ( Y^s ) = 2 \xi (s ) ~~~~~~~~~~(s \in \complex) .
\en 
where $E$ is the expectation operator.
Many other constructions of random variables with the same distribution 
as $Y$  have since been discovered,
involving functionals of the path of a Brownian motion or 
Brownian bridge in $\reals^d$ for $d = 1,2,3$ or $4$.

Our main purpose in this paper is to review this circle of ideas, with 
emphasis on the probabilistic interpretations 
such as \re{bmm}-\re{y1} of various functions which play an important
role in analytic number theory.
For the most part this is a survey of known results,
but the result of Section \ref{renorm} may be new.

Section \ref{analb} reviews the classical analysis underlying \re{y1}, and
offers different analytic characterizations of the probability distribution of
$Y$.  
Section \ref{infdiv} presents various formulae related to the
distributions of the random variables $\sigh$  and $\sighp$ 
defined by 
\eq
\lb{sighdef}
\sigh:= {2 \over \pi^2 } \sum_{n = 1}^\infty { \gamhn \over n^2}
~~~\mbox{ and } ~~~
\sighp:= {2 \over \pi^2 } \sum_{n = 1}^\infty { \gamhn \over (n - \hf)^2},
\en
for independent random variables $\gamhn$ with the 
{\em gamma$(h)$ density}
\eq
\lb{gamhdef}
P(\gamhn \in dx)/dx = \Gamma(h)^{-1} x^{h-1} e^{-x}  ~~~~~~~~~~~~~~(t >0).
\en
Our motivation to study these laws stems from their 
 close connection to the classical functions 
of analytic number theory, and
 their repeated appearances in the study of 
 Brownian motion, which we recall in  Section \ref{brown}.
 For example, to make the connection with the beginning of this introduction,
 one has 
\eq
\lb{sigy}
\sigt\ed {2\over \pi}Y^2
\en
where $\ed$ means equality in distribution. As we discuss in Section 
\ref{brown}, Brownian paths possess a number 
 distributional symmetries, which explain some of the 
 remarkable coincidences in distribution 
implied by the repeated appearances 
of the laws of $\sigh$ and $\sighp$ for various $h$. 
 Section \ref{renorm} shows how one of the probabilistic results of 
 Section \ref{infdiv} leads us to an approximation of the zeta function, 
 valid in
 the entire complex plane, which is similar to an approximation obtained by 
 Sondow \cite{sondow94}. 
 We conclude in Section \ref{sfc} with some 
 consideration of the Hurwitz zeta function and  Dirichlet $L$-functions,
 and some references to other work relating the Riemann zeta function to 
 probability theory. 

\section{Probabilistic interpretations of some classical analytic formulae}
\label{analb}

\subsection{Some classical analysis}
\label{classical}

Let us start with Jacobi's theta function identity
\eq
\lb{jac0}
{1 \over \sqrt{ \pi t}} \sum_{n = -\infty}^\infty e^{- ( n+x)^2/t}
= 
\sum_{n = -\infty}^\infty \cos( 2 n \pi x ) \,\,  e^{- n^2 \pi^2 t }
~~~~~~~~(x\in\real,\, t> 0)
\en
 which is a well known instance of the Poisson summation
formula \cite{bel61}.
This identity equates two different 
expressions for $p_{t/2}(0,x)$,
where $p_{t}(0,x)$ is the fundamental solution of the heat equation on a circle 
identified with $[0,1]$, with initial condition $\delta_0$,
the delta function at zero.
In probabilistic terms, $p_{t}(0,x)$ is the probability density 
at $x \in [0,1]$ of the position of a Brownian motion on the circle
started at $0$ at time $0$ and run for time $t$. 
The left hand expression is obtained
by wrapping the Gaussian solution on the line, while the right hand expression
is obtained by Fourier analysis. 
In particular, \re{jac0} for $x=0$ can be written
\eq
\lb{jac2}
\sqrt{t} \, \theta(t) = \theta(t^{-1}) ~~~~~~~~~~~~~~~~(t >0)
\en
where $\theta$ is the Jacobi theta function
\eq
\lb{jac}
\theta(t):= \sum_{n = -\infty}^\infty \exp( -  n^2 \pi t ) 
~~~~~~~~~~~~~~~~~~(t >0).
\en

For the function $\xi$ defined by \re{xifn},
Riemann  obtained the integral representation
\eq
\lb{riem0}
{4 \xi (  s ) \over s ( s - 1 ) } = \int_0^\infty t^{ {s \over 2} -1} (\theta(t) -1)
dt ~~~~~~~~~~(\Re( s) >  1)
\en
by switching the order of summation and integration and using 
\eq
\lb{gams}
\Gamma(s) = \int_0^\infty x^{s-1} e^{-x } \, dx ~~~~~(\Re s  > 0).
\en
He then deduced his functional equation $\xi(s) = \xi(1-s)$ from \re{riem0} 
and Jacobi's functional equation \re{jac2}.
Following the notation of Edwards \cite[\S 10.3]{edwards74}, 
let 
\eq
\lb{gydef}
G(y):= \theta(y^2) = \sum_{n = -\infty}^\infty \exp( - \pi n^2 y ^2) 
\en
so Jacobi's functional equation
\re{jac2} acquires the simpler form
\eq
\lb{geq}
y G(y) = G(y^{-1}) ~~~~~~~~~~~~~~(y >0).
\en
The function
\eq
\lb{hdef1}
H(y) := 
{ d \over dy } \left[ y^2 {d \over dy } G(y) \right] = 
2 y G'(y) + y^2 G''(y) %%%\mbox{ where }  G(y):= \theta(y^2).
\en
that is
\eq
\lb{hdef}
H(y)= 
4 y^2 \sum_{n=1}^\infty ( 2 \pi^2 n^4 y^2 - 3 \pi n^2 ) e^{- \pi n^2 y^2 } 
\en
satisfies the same functional equation as $G$:
\eq
\lb{heq}
y H(y) = H(y^{-1}) ~~~~~~~~~~~~~~(y >0).
\en
As indicated by Riemann, this allows \re{riem0} to be transformed 
by integration by parts for $\Re s >1$ to yield 
\eq
\lb{edrep}
2 \xi (  s ) = \int_0^\infty y^{s-1} H(y) dy .
\en

It follows immediately by analytic continuation that \re{edrep} serves 
to define an entire function $\xi(s)$ which satisfies Riemann's functional equation $\xi(s) = \xi(1-s)$ for all complex $s$.
Conversely, the functional equation \re{heq} for $H$ is recovered from
Riemann's functional equation for $\xi$ by uniqueness of Mellin transforms.
The representation of $\xi$ as a Mellin transform was used by Hardy to
prove that an infinity of zeros of $\zeta$ lie on the critical line. It is also
essential in the work of \Polya\ \cite{polya74} and Newman \cite{newman76}
on the Riemann hypothesis. But the probabilistic interpretations of
\re{heq} which we discuss in this paper do not appear in these works.

\subsection{Probabilistic interpretation of $2 \xi(s)$}
\label{probint}
As observed by Chung \cite{chu76} and Newman \cite{newman76}, 
$H(y) >0 $ for all $y >0$ 
(obviously for $y \ge 1$, hence too for $y <1$ by \re{heq}).
By \re{xifn}, \re{xieq} and $\zeta(s) \sim (s-1)^{-1}$ as $s \te 1$,
$$ 
2 \xi( 0 ) = 2 \xi (1) = 1,
$$
so formula \re{edrep} for $s=0$ and $s=1$ implies
$$
\int_0^\infty y^{-1} H(y) dy = \int_0^\infty H(y) dy = 1 .
$$
That is to say, the function $y^{-1} H(y)$ is the 
density function of a probability distribution on $(0,\infty)$ with mean 1.
Note that the functional equation \re{heq} for $H$
can be expressed as follows in terms of a random variable $Y$ with this
distribution: for every non-negative measurable function $g$
\eq
\lb{gheq}
E [ g(1/Y) ] = E [ Y g(Y) ].
\en
The distribution of $1/Y$ is therefore identical to the 
{\em size-biased} distribution derived from $Y$. See 
Smith-Diaconis \cite{sd88} for further interpretations of this relation.
The next lemma, which follows from the preceding discussion
and formulas tabulated in Section
\ref{infdiv}, gathers different characterizations
of a random variable $Y$ with this density. Here $Y$ is
assumed to be defined on some probability space $(\Omega, \FF, P)$, with 
expectation operator $E$.

\begin{prp}{lmmy} {\em (\cite{chu76,by87})}
For a non-negative random variable $Y$, each of the following conditions 
{\em (i) - (iv)}
is equivalent to $Y$ having density 
$y^{-1}H(y)$ for $y >0$:

\noindent
{\rm (i)}
\eq
\lb{y1y}
E ( Y^s ) = 2 \xi (s ) ~~~~~~~~~~(s \in \complex);
\en
%\eq
%\lb{y2}
%E [ f(Y) ] = \int_0^\infty f(y) y^{-1} H(y) dy
%\en
%for all non-negative measurable functions $f$;
{\rm (ii)} for $y >0$
\eq
\lb{fdef}
P( Y \le y ) =  G(y) + y G'(y) = -y^{-2} G'(y^{-1})  , 
\en
that is
\eq
\lb{fdef1}
P( Y \le y ) 
= \sum_{n = -\infty}^\infty ( 1 - 2 \pi n^2 y^2 ) e^{- \pi n^2 y^2} 
=  4 \pi y^{-3} \sum_{n=1}^\infty n^2 e^{- \pi n^2/ y^2} ;
\en
{\em (iii)} with $\sigt$ defined by {\em \re{sighdef}}
\eq
\lb{ysigt}
Y \ed \sqrt{ \mbox{$\pi \over 2$} \sigt } ;
\en
{\em (iv)}
\eq
E \left[e^{- \la Y^2} \right] = 
 \left({ \sqrt{\pi \la } \over \sinh \sqrt{\pi \la } }  \right)^2 .
\en
\end{prp}

\section{Two infinitely divisible families}
\label{infdiv}

\noindent
%\vspace*{-6in}
\begin{table}[htbp]
%\noindent\hspace*{-\lmargin}{\hfill\sectit Continuous Distributions \hfill}
%\hfill
\noindent
\hspace*{-1in}\hfill
\begin{tabular}{|c|c|c|}
\hline
\multicolumn{3}{|c|}{}\\ 
\multicolumn{3}{|c|}{\bf Table \tabo\ }\\ 
\multicolumn{3}{|c|}{}\\ 
\hline
& &  \\ 
$\Sigma$&
$\displaystyle{  \sigo:= {2 \over \pi^2} \sum_{n = 1}^\infty { \eps_n \over n^2} }$ &
$\displaystyle{  \sigt:= {2 \over \pi^2} \sum_{n = 1}^\infty { \eps_n + \hat{\eps}_n \over n^2} }$ \\
&  &   \\  \hline
&  &   \\ 
$E [ e^{- \la \Sigma}]$&
$\displaystyle{ {\sqrt{2 \la } \over \sinh \sqrt{2 \la} }}$ &
$\displaystyle{ \left({\sqrt{2 \la } \over \sinh \sqrt{2 \la} }\right)^2}$ \\
&  &   \\  \hline
& &  \\ 
L\'evy density&
$\displaystyle{  \rho_1(x):= {1 \over x} \sum_{n=1}^\infty e^{-\pi^2 n^2 x/2} } $&
$\displaystyle{  2 \rho_1(x)} $ \\
 & &   \\  
\hline
 & &   \\ 
$\displaystyle{f(x):= {d \over dx} P(\Sigma \le x)}$&
$\displaystyle{{d \over dx} \sum_{n = -\infty}^\infty (-1)^n e^{ - n^2 \pi^2 x /2}}$ &
$\displaystyle{{d \over dx} \sum_{n = -\infty}^\infty (1- n^2 \pi^2 x)e^{ - n^2 \pi^2 x /2}}$ \\
&  &   \\  
\hline
& &   \\ 
reciprocal relations&
$\displaystyle{f_1(x) = \left( 2 \over \pi x^3 \right)^{1/2} f_2^\#\left( 4 \over \pi^2 x \right) }$ &
$\displaystyle{f_2(x) = \left( 2 \over \pi x \right)^{5/2} f_2 \left( 4 \over \pi^2 x \right) }$ \\ 
& &   \\  
\hline
& &   \\ 
$\displaystyle{E \left[ g \left( 4 \over \pi^2 \Sigma \right)\right]}$&
$\displaystyle{ \sqrt{\pi \over 2} E \left[ (\sigtp)^{- {1 /2} }  g(\sigtp)\right]}$ &
$\displaystyle{ \sqrt{\pi \over 2} E \left[ (\sigt)^{{1 /2} }  g(\sigt)\right]}$ \\
& &   \\  
\hline
& &   \\ 
$E[\Sigma^s]$&
$\displaystyle{\left( 2^{1-2s} - 1 \over 1 - 2s \right) \left( 2 \over \pi \right)^s 2 \xi (2 s )}$ &
$\displaystyle{\left( 2 \over \pi \right)^s 2 \xi (2 s )}$ \\ 
& &   \\  
\hline
& &   \\ 
$E[\Sigma^n]$&
$\displaystyle{{n! \over (2n)!} ( 2 ^{3n} - 2 ^{n+1} ) (-1)^{n+1} B_{2n}}$ &
$\displaystyle{{n! \over (2n)!} (2n -1) 2^{3n} (-1)^{n+1} B_{2n}}$ \\
& &   \\  
\hline
%& & \\ 
%$xxxx$&
%$\displaystyle{xxx}$ &
%$\displaystyle{xxx}$ \\ 
%& &   \\  
%\hline
\end{tabular}
\label{table1}
\end{table}

\noindent
\begin{table}[htbp]
\noindent
\hspace*{-1in}\hfill
\begin{tabular}{|c|c|}
\hline
\multicolumn{2}{|c|}{}\\ 
\multicolumn{2}{|c|}{\bf Table \tabo\ continued}\\ 
\multicolumn{2}{|c|}{}\\ 
\hline
&   \\ 
$\displaystyle{  \sigop:= {2 \over \pi^2} \sum_{n = 1}^\infty { \eps_n \over (n-\hf)^2} }$ &
$\displaystyle{  \sigtp:= {2 \over \pi^2} \sum_{n = 1}^\infty { \eps_n + \hat{\eps}_n \over ( n - \hf )^2} }$ \\ 
&  \\  \hline
&   \\ 
$\displaystyle{ { 1 \over \cosh \sqrt{2 \la} }}$ &
$\displaystyle{ \left({ 1 \over \cosh \sqrt{2 \la} }\right)^2}$ \\ 
&   \\  \hline
&  \\ 
$\displaystyle{  \rho_1^\# (x):= {1 \over x} \sum_{n=1}^\infty e^{-\pi^2 (n - 1/2)^2 x/2} } $&
$\displaystyle{  2 \rho_1^\#(x)} $ \\
&   \\  
\hline
&   \\ 
$\displaystyle{\pi \sum_{n=0}^\infty (-1)^n ( n + \hf) e^{- (n + {1 \over 2})^2 \pi^2 x/2}}$ &
$\displaystyle{{1 \over 2} \sum_{n = -\infty}^\infty \left( (n + \hf)^2 \pi^2 x - 1 \right) e^{- (n + {1 \over 2})^2 \pi^2 x /2 }}$ \\
&   \\  
\hline
&   \\ 
$\displaystyle{f_1^\#(x) = \left({2 \over \pi x } \right)^{3/2} f_1^\# \left( 4 \over \pi^2 x \right) }$ &
$\displaystyle{f_2^\#(x) = {2 \over \pi} \left( 2 \over \pi x \right)^{3/2} f_1 \left( 4 \over \pi^2 x \right) }$ \\
&   \\  
\hline
&   \\ 
$\displaystyle{ \sqrt{2 \over \pi} E \left[ (\sigop)^{{- 1 /2} }  g(\sigop)\right]}$ &
$\displaystyle{ \left( 2 \over \pi \right)^{3/2} E \left[ (\sigo)^{{- 1 /2} }  g(\sigo)\right]}$ \\
&   \\  
\hline
&   \\ 
$\displaystyle{ \Gamma(s+1) 2^{s+1} \left( 2 \over \pi \right)^{2s+1}L_{ \chi_4}(2 s + 1)}$ &
$\displaystyle{ { (2 ^{2(s+1)} - 1 ) \over s+ 1 } \left( 2 \over \pi \right)^{s+1} \xi ( 2(s+1) ) }$ \\
&   \\  
\hline
&   \\ 
$\displaystyle{{n! \over (2n)! } 2^n (-1)^n E_{2n}}$ & 
$\displaystyle{{ ( 2^{2n+2} - 1) 2^{3 n + 1 } \, n! \over (n+1) (2n)! } (-1)^n B_{2n+2}} $\\
&   \\  
\hline
%&   \\ 
%$\displaystyle{xxx}$ &
%$\displaystyle{xxx}$ \\
%&   \\  
%\hline
\end{tabular}
\label{table2r}
\end{table}
This section presents an array of results regarding the probability
laws on $(0,\infty)$ of the random variables $\sigh$ and $\sighp$ defined by
\re{sighdef}
with special emphasis on results for $h = 1$ and $h = 2$, which
are summarized by Table \tabo.
Each column of the table presents features of the law of one of the
four sums $\Sigma = \sigo, \sigt, \sigop$ of $\sigtp$. 
Those in the $\sigt$ column can be read from Proposition \ref{lmmy}, 
while the formulae in other columns provide to analogous results
for $\sigo$, $\sigop$ and $\sigtp$ instead of $\sigt$.
While the Mellin transforms of $\sigo$, $\sigt$ and
$\sigtp$ all involve the function $\xi$ associated with the
Riemann zeta function, the Mellin transform of $\sigop$ involves instead
the Dirichlet $L$-function associated with the quadratic character modulo 4,
that is 
\eq
\lb{oddser}
L_{\chi_4}(s):= \sum_{n=0}^\infty { (-1)^n \over ( 2 n + 1 ) ^s }  
%= 4^{-s} [\zeta(s, \qr) - \zeta(s, \tqr)]  
~~~~(\Re s > 0 ).
\en
We now discuss the entries of Table \tabo\ row by row.

\subsection{Laplace transforms and \Lev\ densities}
Recall from \re{sighdef} that
$\sigh:= { 2 \over \pi^2 } \sum_{n = 1}^\infty {\gamhn \over n^2}$
where the distribution of the independent gamma$(h)$ variables 
$\gamhn$ is characterized by the Laplace transform
\eq
\lb{gamhlt}
E[ \exp ( - \la \gamhn ) ] = (1 + \la )^{-h} .
\en
Euler's formulae 
\eq
\lb{esinh}
\sinh z = z \prod_{n=1}^\infty \left( 1 + { z^2 \over { n^2 \pi ^2 } }\right)
~~~~\hbox{and}~~~~~~
\cosh z = \prod_{n=1}^\infty \left( 1 + {  z^2 \over { (n - \hf)^2 \pi ^2 } }\right)
\en
 allow the following evaluation 
\cite{watson61,chu82}:
for $\Re (\la ) > 0$    %%%- \pi^2/2$
\eq
\lb{y44}
E \left[e^{- \la \sigh } \right] = E \left[ 
\prod_{n=1}^\infty \exp \left( - {2 \la \gamhn\over n^2 \pi^2 }
 \right) \right]
= \prod_{n=1}^\infty 
\left( 1 +   { 2 \la \over n^2 \pi^2} \right)^{-h} 
= \left({ \sqrt{2 \la } \over \sinh \sqrt{2 \la } }  \right)^h 
\en
and similarly
\eq
\lb{coshlt}
E  [ \exp ( - \la \sighp ) ] = \left({1 \over \cosh \sqrt{2 \la } }\right)^h .
\en

\noindent{\bf \Lev\ densities $\rho(x)$.}
A probability distribution $F$ on the line is called {\em infinitely
divisible} if for each $n$ there exist independent random variables
$T_{n,1} , \ldots, T_{n,n} $ with the same distribution such that
$\sum_{i = 1}^n T_{n,i}$ has distribution $F$. 
According to the {\em \LK\ representation}, which has a well known
interpretation in terms of Poisson processes \cite[\S II.64]{wil79}, 
a distribution $F$ concentrated on $[0,\infty)$ is infinitely
divisible if and only if
its Laplace transform $\lt(\la):= \int_0^\infty e^{-\la t} F (dt )$
admits the representation
\eq
\lb{lkform}
\lt(\la) = \exp \left( -c \la   - \int_0^\infty (1 - e^{-\la x} ) \nu (dx) \right)
~~~~~~~(\la \ge 0 )
\en
for some $c \ge 0$ and some positive 
measure $\nu$ on $(0,\infty)$, called the 
{\em \Lev\ measure} of $F$, with $c$ and $\nu$ uniquely determined by $F$.

If $\nu(dx) = \rho(x) dx$ 
%%for some non-negative function $\rho$, 
then $\rho(x)$ is called the {\em \Lev\ density} of $F$.  
It is elementary that for $\gamh$ with gamma$(h)$ distribution and $a >0$ 
the distribution of $\gamh/a$ is infinitely divisible with
\Lev\ density $hx^{-1} e^{-ax}$. It follows easily that for $a_n >0$ with
$\sum_n a_n < \infty$ and independent gamma$(h)$ variables
$\gamhn$ the distribution of $\sum_n \gamhn /a_n$ is
infinitely divisible with \Lev\ density $h x^{-1} \sum_n e^{-a_nx}$. Thus
for each $h >0$ the laws of $\sigh$ and $\sighp$ are infinitely divisible,
with the \Lev\ densities indicated in the Table for $h = 1$ and $h =2$.

We note that Riemann's formula \re{riem0} for $s = 2p$
can be interpreted as 
an expression for the $p$th moment $\mu_1(p)$ of the \Lev\ density $\rho_1$ 
of $\sigo$: for $\Re p > \hf$
\eq
\lb{levzeta}
\mu_1(p):=
\int_0^\infty t^{p} \rho_1(t) dt = 
\left( 2 \over \pi \right)^{p} {1 \over 2} \int_0^\infty y^{{p } - 1 }( \theta(y) - 1 ) dy =
{2^p \over \pi^{2p}} \Gamma(p) \zeta(2p)  .
\en

\subsection{Probability densities and reciprocal relations}

By application of the negative binomial expansion
\eq
\lb{negbin}
{1 \over (1-x)^{h}} = {1 \over \Gamma(h)} \sum_{n=0}^\infty { \Gamma(n+h) \over \Gamma(n+1)} x^n  ~~~~~~~~(h >0, |x| < 1 )
\en
there is the expansion
\eq
\lb{sinhexp}
\left( t \over \sinh t   \right)^h
=  {2^h t^h e^{-t h } \over ( 1 - e^{-2t})^{h}}
= { 2^h t^h \over \Gamma(h) }
\sum_{n=0}^\infty { \Gamma (n+h) \over \Gamma(n+1) } e^{ - (2n + h ) t }
\en
which corrects two typographical errors in \cite[(3.v)]{by87}, and
\eq
\lb{coshexp}
\left( 1 \over \cosh t  \right)^h
=  {2^h e^{-t h } \over ( 1 + e^{-2t})^{h}}
= { 2^h \over \Gamma(h) } \sum_{n=0}^\infty (-1)^n { \Gamma (n+h) \over \Gamma(n+1) } e^{- (2n + h ) t } .
\en
The Laplace transform \re{coshlt} can be inverted by applying
the expansion \re{coshexp} and
inverting term by term using \Lev's formula \cite{lev39} 
%for the Laplace density of the one-sided stable distribution of index 1/2
\eq
\lb{lev}
\int_0^\infty { a \over \sqrt{2 \pi t^3} } 
e^{-a^2 / (2 t)} e^{-\la t } dt = e^{- a \sqrt{2 \la }} .
\en
Thus there is the following expression for the density 
$f_h^\#(t):= P(\sighp \in dt)/dt$: for arbitrary real $h >0$:
\eq
\lb{cosh3x}
f_h^\# (t) = {2^h \over \Gamma(h) } \sum_{n = 0}^\infty (-1)^n 
{ \Gamma ( n + h ) \over \Gamma(n+1) } {(2n+h) \over \sqrt{2 \pi t^3}}
\exp \left( - {(2 n + h )^2 \over 2t } \right) .
\en
A more complicated formula for $f_h(t)$ was obtained from
\re{sinhexp} by the same method in \cite[(3.x)]{by87}.
The formulae for the densities of $\sigh$ and $\sighp$ displayed
in Table \tabo\ for $h = 1$ and $h = 2$ can be obtained 
using the reciprocal relations of Row 5. The self-reciprocal relation
involving $\sigt$ is a variant of \re{heq}, while that 
involving $\sigop$, which was observed by Ciesielski-Taylor \cite{ct62},
is an instance of another application of the Poisson
summation formula which is recalled as \re{fethetachi} in Section \ref{hur}.
Lastly, the reciprocal relation involving the
densities $f_1$ of $\sigo$ and $f_2^\#$ of $\sigtp$ amounts
to the identity
\eq
\lb{kolcdf}
P( \sigo \le x ) =  \sum_{n = -\infty}^\infty (-1)^n e^{ - n^2 \pi^2 x /2}
= 
{ \sqrt{2 \over \pi x} }\sum_{n= -\infty}^{\infty}e^{-2 (n+ \hf)^2/x} 
\en
where the second equality is read from \re{jac0} with 
$x = 1/2$ and $t$ replaced by $x/2$.

\noindent{\bf Formulae for $E\left[g\left(4 \over \pi^2 \Sigma \right) \right]$.}
 These formulae, valid for an arbitrary non-negative Borel 
 function function $g$, are integrated forms of the reciprocal relations,
 similar to \re{gheq}.

\subsection{Moments and Mellin transforms}
It is easily shown that the distributions of $\sigh$ and $\sighp$ have 
moments of all orders (see e.g.  Lemma \ref{conv}).
The formulae for the Mellin transforms $E(\Sigma^s)$
can all be obtained by term-by-term integration
of the densities for suitable $s$, followed by analytic continuation.
According to the self-reciprocal relation  for $\sigop$,
for all $s \in \complex$
\eq
\lb{msigop}
E \left[ \left(\mbox{ $\pi \over 2 $} \sigop \right) ^{s} \right]
= E \left[ \left(\mbox{ $\pi \over 2 $} \sigop \right) ^{- {1 \over 2} - s} \right] .
\en
Using the formula for $E((\sigop)^s)$ in terms of
$L_{\chi_4}$ defined by \re{oddser}, given in Table 1,
%\eq
%\lb{msigopp}
%E \left[ \left(\mbox{ $\pi \over 2 $} \sigop \right) ^{s} \right]
%= \left(2 \over \pi \right)^{s +1 } \Gamma(s +1) 2^{s+1} L_{\chi_4} (2 s +1)
%\en
%first for $\Re s > - \hf$, then for all $s \in \complex$ by analytic
%continuation.
we see that if we define
\eq
\lb{chitil}
\Lambda_{\chi_4} (t):=
 E \left[ \left(\mbox{ $\pi \over 2$ } 
 \sigop \right) ^{{t-1 \over 2}} \right] = 
\Gamma \left( {t +1 \over 2 } \right)
 \left( 4 \over \pi \right)^{{t+1 \over 2 }} L_{\chi_4}(t) 
\en
then \re{msigop} amounts to the functional equation
\eq
\lb{chieq}
\Lambda_{\chi_4} (t) = \Lambda_{\chi_4} (1-t) ~~~~~~~~~~(t \in \complex).
\en
This is an instance of the general functional equation 
for a Dirichlet $L-$function, which is recalled as \re{Lfe} in 
Section \ref{sfc}.

\noindent{\bf Positive integer moments $E(\Sigma^n)$.}
These formulas are particularizations of the preceding Row, using 
 the classical evaluation  of $\zeta(2n)$ in terms of the Bernoulli numbers
 $B_{2n}$.
The result for $\sigop$ involves the
{\em Euler numbers} $E_{2n}$ defined by the expansion
\eq
\lb{eulern}
{1 \over \cosh(z)} = { 2 \over e^{z} + e^{-z} } = 
\sum_{n=0}^\infty { E_{2n} } { z^{2n} \over (2n)! }
\en

\noindent{\bf A  multiplicative relation.}
Table \tabo\ reveals the following remarkably simple relation:
\eq
\lb{sigot}
E(\sigo ^s ) =  \left({ 2 ^{1- 2s} - 1 \over 1 - 2 s }  \right) E( \sigt ^s ) 
\en
where the first factor on the right side is evaluated by continuity for
$s = 1/2$.
By elementary integration, this factor can be interpreted as follows:
\eq
\lb{wws}
\left({ 2 ^{1- 2s} - 1 \over 1 - 2 s }  \right)= E ( W^{-2s} )
\en
for a random variable $W$ with uniform distribution on $[1,2]$. 
Thus \re{sigot} amounts to the
the following identity in distribution:
\eq
\lb{prodid}
\sigo \ed W^{-2} \sigt 
\en
where $W$ is assumed independent of $\sigt$.
An equivalent of \re{prodid} was interpreted in terms of Brownian
motion in \cite[(4)]{py97kt}.

Note that \re{prodid} could be rewritten as
$\sigo \ed  W^{-2} (\sigo + \sigoh)$ for
$\sigo$, $\sigoh$ and $W$  independent random variables, with
$\sigoh$ having the same distribution as $\sigo$, and $W$ uniform on $[1,2]$.
By consideration of positive integer moments, this property uniquely
characterizes the distribution of $\sigo$ among all distributions
with mean $1/3$ and finite moments of all orders.

\subsection{Characterizations of the distributions of $\sigt$ and $\sigtp$}
\label{charact}
As just indicated, the identity \re{prodid} allows a simple
probabilistic
characterization of the distribution of $\sigo$. The following Proposition
offers similar characterizations of the distributions of $\sigt$ and $\sigtp$.
\begin{prp}{pccosh}
Let $X$ be a non-negative random variable, and let 
$\tX$ denote a random variable such that
$$P(\tX \in dx ) = x P( X \in dx )/E(X).$$
\noindent{\em (i)}  $X$ is distributed as $\sigt$ if and only if
$E(X) = 2/3$ and 
\eq
\lb{xhx}
\tX \ed X + \HH \tX
\en
where $X, \HH$ and $\tX$ are independent, with
$$
P(\HH \in dh ) = ( h^{-1/2} - 1) dh ~~~~~~~~~~~(0 < h < 1)
$$
\noindent{\em (ii)}  $X$ is distributed as $\sigtp$ if and only if
$E(X) = 2$ and 
\eq
\lb{xcosh}
\tX \ed X + U^2 \hat{X}
\en
where $X$, $U$ and $\hat{X}$ are independent,
with $\hat{X}$ distributed as $X$ and $U$ uniform on $[0,1]$.
\end{prp}
For the proof of this Proposition, note  that \re{xhx} (or \re{xcosh})
 imply that
the Laplace transform of
$X$ satisfies an integro-differential equation,
 whose only solution is given by the appropriate function.
 The ``only if'' part of (i) appears in
\cite[p. 26]{yor97z}. Details of the remaining parts will be provided 
elsewhere.
As remarked by van Harn and Steutel \cite{vhs95}, 
it is an easy consequence of the \LK\ formula that for non-negative 
random variables $X$ and $Y$ the equation $X^* \ed X + Y$ is satisfied
for some $Y$ independent of $X$ if and only if the law of $X$ is infinitely
divisible. As discussed in \cite{vhs95,arrgold98},
the distribution of $X^*$, known as the {\em size-biased}
or {\em length-biased} distribution of $X$, has a natural interpretation
in renewal theory. 

\section{Brownian interpretations}
\label{brown}

It is a remarkable fact that the four distributions considered in Section
\ref{infdiv} appear in many different problems concerning Brownian motion and
related stochastic processes. These appearances are partially
explained by the relation of these distributions to Jacobi's theta function, 
which provides a solution to the heat equation \cite{bel61,ehren87}, 
and is therefore related to Brownian motion \cite[\S 5.4]{kni88}.
We start by introducing some basic notation for Brownian motion and
Bessel processes, then present the main results in the form of
another table.

\subsection{Introduction and Notation}
\label{oned}

Let $\bmbeta:= (\bmbeta_t, t \ge 0 )$ be a standard
one-dimensional Brownian motion, that is a stochastic process with continuous
sample paths and independent increments such that $\bmbeta_0 = 0$, and for 
all $s,t >0 $ the random variable $\bmbeta_{s+t} - \bmbeta_s$ has a Gaussian distribution with mean $E( \bmbeta_{s+t} - \bmbeta_s ) = 0$ and 
mean square $E[( \bmbeta_{s+t} - \bmbeta_s )^2] = t$, meaning that for all real $x$
$$
P( \bmbeta_{s+t} - \bmbeta_s \le x ) = 
{1 \over \sqrt{2 \pi t } } \int_{-\infty}^x e^{-y^2/(2t) }\,dy.
$$
Among continuous time stochastic processes, such as semimartingales,
processes with independent increments, and Markov processes,
Brownian motion is the paradigm of a stochastic process with continuous
paths.
In particular, among processes with
stationary independent increments, the Brownian motions
$(\sigma B_t +\mu t, t \ge 0 )$ for $\sigma >0, \mu \in \reals$
are the only ones with almost surely continuous paths \cite[I.28.12]{rw94}.
Brownian motion arises naturally as the limit in distribution as $n \te \infty$
of a rescaled random walk process $(S_n, n = 0,1, \ldots)$ where
$S_n = X_1 + \cdots + X_n$ for independent random variables $X_i$ with
some common distribution with mean $E(X_i) = 0$ and variance $E(X_i ^2 ) = 1$.
To be more precise, let the value of $S_r$ be extended to all real $r \ge 0$
by linear interpolation between integers. 
With this definition of $(S_r, r \ge 0)$ as a random continuous function,
it is known that no matter what the distribution of the $X_i$
with mean $0$ and variance $1$, as $n \te \infty$
\eq
\lb{donsk}
\left( {S_{nt} \over \sqrt{n}}, t \ge 0 \right) \convd (\bmbeta_t, t \ge 0)
\en
in the sense of weak convergence of probability distributions on the
path space $C[0,\infty)$. In particular, convergence of
finite dimensional distributions in \re{donsk} follows easily from
the central limit theorem, which is the statement of convergence of 
one dimensional distributions in \re{donsk},
that is for each fixed $t>0$
\eq
{S_{nt} \over \sqrt{n}} \convd \bmbeta_t
\en
where $\convd$ denotes weak convergence of probability distributions on
the line. Recall that,
 for random variables $W_n, n = 1,2, \ldots$ and $W$ such that
$W$ has a continuous distribution function $x \mapsto P(W \le x)$,
$W_n \convd W$ means $P(W_n \le x ) \te P(W\le x)$ for all real $x$.
See \cite{bil68,ry99} for background.

Let $(b_t, 0 \le t \le 1)$ be a {\em standard Brownian bridge}, that is the
centered Gaussian process with the conditional distribution 
of $(\bmbeta_t, 0 \le t \le 1)$ given
$\bmbeta_1 = 0$. Some well known alternative descriptions of the distribution
of $b$ are \cite[Ch. III, Ex (3.10)]{ry99} 
\eq
(b_t , 0 \le t \le 1) \ed  ( \bmbeta_t - t \bmbeta_1, 0 \le t \le 1)
\ed ((1 - t ) \bmbeta_{t/(1-t)}, 0 \le t \le 1)
\en
where $\ed$ denotes equality of distributions on the path space $C[0,1]$,
and the rightmost process is defined to be $0$ for $t = 1$.
According to a fundamental result in the theory of
non-parametric statistics \cite{doob49,sw86}, 
the Brownian bridge arises 
in another way from the asymptotic behaviour of the 
{\em empirical distribution}
$$
F_n(x):= {1 \over n } \sum_{ k = 1}^n 1( X_k \le x )
$$
where the $X_k$ are now supposed independent with common distribution
$P(X_i \le x ) = F(x)$ for an arbitrary continuous distribution function $F$.
As shown by Kolmogorov \cite{kol33}, 
the
distribution of $\sup_x |F_n(x) - F(x)|$ is the same 
no matter what the choice of $F$, and for all real $y$
\eq
\lb{kolm}
\lim_{n \te \infty } P( \sqrt{n } \sup_x |F_n(x) - F(x)| \le y )
= \sum_{n = -\infty}^\infty (-1)^n e^{- 2 n^2 y^2 } 
\en
For $F$ the uniform distribution on $[0,1]$, 
so $F(t) = t$ for $0 \le t \le 1$, it is known that 
\eq
\lb{edf} 
(\sqrt{n } (F_n(t) - t), 0 \le t \le 1 ) \convd (b_t, 0 \le t \le 1) .  
\en 
As a well known consequence of \re{edf}, Kolmogorov's limiting distribution 
in \re{kolm} is identical to the distribution of $\max_{0\le t \le 1}|b_t|$. 
On the other hand, as observed by Watson \cite{watson61}, Kolmogorov's 
limit  distribution function in \re{kolm} is identical to that of
${\pi \over 2 } \sqrt{ \TT_1}$.
Thus we find the first appearance of the law of $\TT_1$ as the
law of a functional of Brownian bridge.

To put this in terms of random walks, if $(S_n)$ is a 
{\em simple random walk},
meaning $P(X_i = +1) = P(X_i =  -1) = 1/2$, then
\eq
\lb{cinv}
\left( \left.{S_{2nt} \over \sqrt{ 2 n}}, 0 \le t \le 1 \right| S_{2n} = 0 \right) \convd (b_t, 0 \le t \le 1)
\en
where on the left side the random walk is conditioned to return to zero
at time $2n$, and on the right side the
Brownian motion is conditioned to return to zero at time 1. 
Thus
\eq
\lb{rw1}
\left( \left. {1 \over \sqrt{ 2 n} } \max_{0 \le k \le 2n} |S_k| \right| S_{2n} = 0 \right) 
\convd \max_{0 \le t \le 1} |b_t| \ed {\pi \over 2 } \sqrt{ \TT_1}
\en
where the equality in distribution summarizes the conclusion of the previous
paragraph.
In the same vein, Gnedenko \cite{gned54} derived 
another asymptotic distribution from random walks, which
can be interpreted in terms of Brownian bridge as
\eq
\lb{gngi1}
\left( \left. {1 \over \sqrt{ 2 n} } \left[\max_{0 \le k \le 2n} S_k - \min_{0 \le k \le 2n} S_k \right] \right| S_{2n} = 0 \right) 
\convd 
\max_{0 \le t \le 1} b_t - \min_{0 \le t \le 1} b_t  \ed { \pi \over 2 } \sqrt{ \TT_2} .
\en
The equalities in distribution in both \re{rw1} and \re{gngi1} can be
deduced from the formula
\eq
\lb{smirdoob}
P(\min_{0 \le u \le 1} b_u \ge - a, \max_{0 \le u \le 1} b_u  \le b) 
=\sum_{k= - \infty}^\infty 
e^{ - 2 k^2 ( a + b )^2 }
- \sum_{k= - \infty}^\infty e^{ - 2 [ b + k ( a + b )]^2 }
\en
of Smirnov \cite{smirn39} and Doob \cite{doob49}.
Kennedy \cite{ken76} found that these distributions appear again
if the random walk is conditioned instead on the event $(R=2n)$ or $(R > 2n)$,
where 
$$R: =\inf\{n\geq 1:  S_n=0\}$$ 
is the time of the first return to zero by the random walk.
Thus
\eq
\lb{exclim}
\left( \left. {1 \over \sqrt{ 2 n} } \max_{0 \le k \le 2n} |S_k |
\right| R = 2n \right) 
\convd 
\max_{0 \le t \le 1} e_t \ed { \pi \over 2 } \sqrt{ \TT_2} 
\en
where $(e_t, 0 \le t \le 1)$ denotes a {\em standard Brownian excursion},
that is the process with continuous sample paths 
defined following \cite{kaigh76,ken76,dim77} by the limit in distribution on $C[0,1]$
\eq
\lb{excinv}
\left( \left.{|S_{2nt}| \over \sqrt{ 2 n}}, 0 \le t \le 1 \right| R = 2n \right) \convd (e_t, 0 \le t \le 1) .
\en
A satisfying explanation of the identity in distribution between
the limit variables featured in \re{gngi1} and \re{exclim} is provided by the
following identity of distributions on $C[0,1]$ due to
Vervaat \cite{ver79}: 
\eq
\lb{verv}
(e_u, 0 \le u \le 1) \ed (b_{\rho + u {\rm (mod 1)}}- b_\rho, 0 \le u \le 1 ) 
\en
where $\rho$ is the almost surely unique
time that the Brownian bridge $b$ attains its minimum value.
As shown by \Tak\cite{tak58} and Smith-Diaconis \cite{sd88}, either of the
approximations \re{gngi1} or \re{exclim} can be used
to establish the differentiated form \re{fdef} of 
Jacobi's functional equation \re{geq} by a discrete 
approximation argument involving quantities of probabilistic interest.
See also \Pol\ \cite{polya27} for a closely related
proof of Jacobi's functional equation 
based on the local normal approximation to the binomial distribution.

In the same vein as \re{gngi1} and \re{exclim} there is the result of
\cite{ken76,dim77} that
\eq
\lb{meanlim}
\left( \left. {1 \over \sqrt{ 2 n} } \max_{0 \le k \le 2n} |S_k |
\right| R > 2n \right) 
\convd 
\max_{0 \le t \le 1} m_t \ed \pi \sqrt{ \TT_1}
\en
where $(m_t, 0 \le t \le 1)$ denotes a {\em standard Brownian meander},
defined by the limit in distribution on $C[0,1]$
\eq
\lb{meaninv}
\left( \left.{|S_{2nt}| \over \sqrt{ 2 n}}, 0 \le t \le 1 \right| R > 2n \right) \convd (m_t, 0 \le t \le 1) .
\en
The
surprising consequence of \re{rw1} and \re{meanlim},
that $\max_{0 \le t \le 1} m_t \ed 2 \max_{0 \le t \le 1} |b_t |$,
was explained in \cite{by87} by a transformation of bridge
$b$ into a process distributed like the meander $m$.
For a review of various transformations relating Brownian
bridge, excursion and the meander see \cite{bp92}.

\subsection{Bessel processes}

The work of %%%Ray \cite{xx}, Knight \cite{xx} and 
Williams \cite{wil70,wil74,wil79} shows how the study of excursions 
of one-dimensional 
Brownian motion leads inevitably to descriptions of these excursions
involving higher dimensional Bessel processes.
For $d = 1,2, \ldots$ let
$R_d:= (R_{d,t}, t \ge 0 )$ be the {\em $d$-dimensional Bessel process $\besd$},
that is the non-negative process defined by the 
radial part of a $d$-dimensional Brownian motion:
$$
R_{d,t}^2 := \sum_{i =1}^d \bmB_{i,t} ^2
$$
where $(\bmB_{i,t}, t \ge 0)$ for $i = 1,2, \ldots$ is
a sequence of independent one-dimensional Brownian motions. 
%$$
%T_d:= \inf \{ t : R_d (t ) = 1 \}
%$$
%denote the hitting time of $1$ by $R_d$.
Note that each of the processes $X = B$, and $X = R_d$ for
any $d \ge 1$, has the {\em Brownian scaling property}:
\eq
\lb{bscale}
(X_{u}, u \ge 0 ) \ed (\sqrt {c} X_{u/c}, u \ge 0 )
\en
for every $c>0$, where $\ed$ denotes equality in distribution
of processes.
For a process $X = (X_t, t \ge 0 )$ let $\overline{X}$ and $\underline{X}$
denote the past maximum and past minimum processes derived from $X$,
that is
$$
\overline{X}_t := \sup_{0 \le s \le t} X_s;~~ \underline{X}_t := \inf_{0 \le s \le t} X_s .
$$
Note that if $X$ has the Brownian scaling property \re{bscale} then 
so too do $\overline{X}, \underline{X}$, and $\overline{X} - \underline{X}$. 
For a suitable process $X$, let
$$
(L_{t}^x(X), t \ge 0, x \in \reals)
$$
be the {\em process of local times of $X$}
defined by the occupation density formula
\eq
\lb{ltdef}
\int_0^t f(X_s) ds = \int_{-\infty}^\infty f(x) L_t^x(X) dx
\en
for all non-negative Borel functions $f$, and 
almost sure joint continuity in $t$ and $x$.
See \cite[Ch. VI]{ry99} for background, and proof of the existence of
such a local time process for $X = B$ and $X = R_d$ for any $d \ge 1$.

Let $r_d := (r_{d,u}, 0 \le u \le 1)$ denote the
{\em $d$-dimensional Bessel bridge} defined by conditioning
$R_{d,u}, 0 \le u \le 1$ on $R_{d,1} = 0$.
Put another way, $r_d^2$ is the sum of squares of $d$ independent copies 
of the standard Brownian bridge.

\subsection{A table of identities in distribution}
\label{table}

We now discuss the meaning of Table \tabt, which presents a number of
known identities in distribution.
The results are collected from the work of numerous authors, including
Gikhman \cite{gik}, Kiefer \cite{kie58}, Chung \cite{chu76},
Biane-Yor \cite{by87}. See also \cite{wzeta,py95agr,py97max}.
In the following sections we review briefly the main arguments underlying the
results presented in the table.

Each column of the table displays a list of random variables 
with the distribution determined by the Laplace transform in Row 0. 
Each variable in the second column is distributed as the sum of 
two independent copies of any variable in the first column, 
and each variable in the fourth column is distributed as 
the sum of two independent copies of any variable in the third column.
The table is organized by rows of variables which are analogous in some 
informal sense. 
The next few paragraphs introduce row by row
the notation used in the table, with pointers to
explanations and attributions in following subsections.
Blank entries in the table mean we do not know any construction of a 
variable with the appropriate distribution which respects the
informal sense of analogy within rows, with the following exceptions.
Entries for Rows 4 and 6 of the 
$\sigt$ column could be filled 
like in Row 3 as the sums of two independent copies of 
variables in the $\sigo$ column of 
the same row, but this would add nothing to the content of the table.
The list of variables involved is by no means exhaustive:
for instance, according to \re{gngi1} the variable
$(4/\pi^2) (\overline{b}_1 - \underline{b}_1)^2$ could be added to the second
column.  Many more constructions are possible involving 
Brownian bridge and excursion, some of which we mention in following 
sections. It is a consequence of its construction, that 
each column of the table exhibits a family of random variables with the same
distribution. Therefore it is a natural problem, coming from 
 the philosophy of
 ``bijective proofs''  in enumerative combinatorics (see e.g.
Stanley \cite{stanleyv2}), to try  giving a direct argument for each
distributional identity, not using the explicit computation of the
distribution. Many such arguments can be given, relying on
distributional symmetries of Brownian paths, or some deeper
 results such as the Ray-Knight theorems.
However, some identities remain for which we do not have any
such argument at hand. As explained in Section \ref{max},
some of these identities are equivalent
to the functional equation for the Jacobi theta (or the Riemann zeta) function.

\noindent {\bf Rows 0 and 1.}
Row 0 displays the Laplace transforms in $\la$ of the distributions
of the variables in Row 1, that is $\sigo$, $\sigt$, 
$\sigop$ and $\sigtp$, as considered in previous sections.

\noindent {\bf Row 2.} 
Section \ref{sqbes} explains why the distributions of the random variables 
$\int_0^1 r_{d,u}^2 du$ and $\int_0^1 R_{d,u}^2 du$ 
for $d = 2$ and $d=4$, are as indicated in this row.

\noindent {\bf Row 3.} 
Most of the results of this row are discussed in Section \ref{fpt}.
Here 
$$
T_a(X):= \inf \{ t : X_t = a \}
$$
is the hitting time of $a$ by a process $X$, and
$\hat{R}_d$ is
%T_1(R_d):= \inf \{ t : R_d (t ) = 1 \}
%$T_1(\hat{R}_d)$ is the
%hitting time of $1$ by $\hat{R}_d$, 
an independent copy of the Bessel process $R_d$.
Note that $R_1 := |B|$ is just Brownian motion with reflection at $0$,
and
$T_1 ( \overline{B} - \underline{B} )$ is the first time that the range of
the Brownian $B$ up to time $t$ is an interval of length $1$. 
The result that
$4 T_1 ( \overline{B} - \underline{B} )$ has Laplace 
transform $1/\cosh ^2 \sqrt{2 \la}$ is due to Imhof \cite{imh85}.
See also Vallois \cite{val92,val95}, Pitman \cite{jp95cyc} and Pitman-Yor \cite{py97kt}
for various refinements of this formula.

\noindent {\bf Rows  4 and 5} 
These rows, which involve the distribution of the maximum of various
processes over $[0,1]$, are discussed in Section \ref{max}.

\noindent {\bf Row 6} 
Here $\overline{m}_1$ is the maximum of the standard Brownian meander
$(m_u, 0 \le u \le 1)$. This entry is read from \re{meanlim}.

\noindent {\bf Row 7.} 
The first two entries are obtained from their relation to the 
first two entries in Row 5, that is the equalities in distribution
$$
\int_0^1 { du \over m_u } \ed 2 \overline{r}_{1,1}
\mbox{  and }
\int_0^1 { du \over r_{3,u} } \ed 2 \overline{r}_{3,1} .
$$
These identities follow from descriptions of the local time processes
$(L_1^{x}(r_{d}), x \ge 0 )$ for $d =1$ and $d =3$,
which involve $m$ for $d=1$ and $r_3$ for $d=3$,
as presented in Biane-Yor \cite[Th. 5.3]{by87}. See also 
\cite[Cor. 16]{jp97sde} for another derivation of these results.
The last two entries may be obtained through their relation to the last two
entries of Row 2. More generally, there is the identity 
$$
{1 \over 2} \int_0^1 { ds \over R_{d,s} } \ed \left( \int_0^1 R_{2d-2,s} ^2 \,ds\right)^{-1/2} ~~~~~~~(d >1)
$$
which can be found in Biane-Yor \cite{by87} and Revuz-Yor 
\cite[Ch. XI, Corollary 1.12 and p. 448]{ry99}.

\noindent {\bf Row 8.} 
Here $\tau_1:= \inf \{t : L_{t}^0(B) = 1 \}$ where $R_1 = |B|$ and
$(L_{t}^x(B), t \ge 0, x \in \reals)$ is the local time process of $B$
defined by \re{ltdef}.
The distribution of 
$\tau_1 / \overline{R}_{1, \tau_1}^2$ was identified with that of $4 T_1(R_3)$
by Knight \cite{kni88}, while the distribution of
$\tau_1 / ( \overline{B}_{\tau_1} - \underline{B}_{\tau_1})^2 $ was identified
with that of $T_1(R_3) + T_1 (\hat{R}_3)$ by Pitman-Yor \cite{py97kt}.
The result in the third column can be read from Hu-Shi-Yor
\cite[p. 188]{hsy96}.

\noindent
%\vspace*{-6in}
\begin{table}[htbp]
%\noindent\hspace*{-\lmargin}{\hfill\sectit Continuous Distributions \hfill}
%\hfill
\noindent
\hspace*{-1in}\hfill
\begin{tabular}{|c|c|c|c|c|}
\hline
\multicolumn{5}{|c|}{}\\ 
\multicolumn{5}{|c|}{\bf Table \tabt\ }\\ 
\multicolumn{5}{|c|}{}\\ 
\hline
& & & &   \\ 
0)&
$\displaystyle{ {\sqrt{2 \la } \over \sinh \sqrt{2 \la} }}$ &
$\displaystyle{ \left({\sqrt{2 \la } \over \sinh \sqrt{2 \la} }\right)^2}$ &
$\displaystyle{ { 1 \over \cosh \sqrt{2 \la} }}$ &
$\displaystyle{ \left({ 1 \over \cosh \sqrt{2 \la} }\right)^2}$ \\ 
& & & &   \\  \hline\hline 
& & & &   \\ 
1)&
$\displaystyle{  \sigo:= {2 \over \pi^2} \sum_{n = 1}^\infty { \eps_n \over n^2} }$ &
$\displaystyle{  \sigt:= {2 \over \pi^2} \sum_{n = 1}^\infty { \eps_n + \hat{\eps}_n \over n^2} }$ &
$\displaystyle{  \sigop:= {2 \over \pi^2} \sum_{n = 1}^\infty { \eps_n \over (n-\hf)^2} }$ &
$\displaystyle{  \sigtp:= {2 \over \pi^2} \sum_{n = 1}^\infty { \eps_n + \hat{\eps}_n \over ( n - \hf )^2} }$ \\ 
& & & &   \\  \hline
& & & &   \\ 
2)&
$\displaystyle{  \int_0^1 r_{2,u}^2 du}$ &
$\displaystyle{  \int_0^1 r_{4,u}^2 du}$ &
$\displaystyle{  \int_0^1 R_{2,u}^2 du}$ &
$\displaystyle{  \int_0^1 R_{4,u}^2 du}$ \\
& & & &   \\  \hline
& & & &   \\ 

3)&
$\displaystyle{  T_1(R_3) }$ &
$\displaystyle{  T_1(R_3) + T_1(\hat{R}_3)}$ &
$\displaystyle{  T_1(R_1) }$ &
$\displaystyle{ 4 T_1(\overline{B} - \underline{B}) }$ \\
& & & &   \\  \hline
& & & &   \\ 

4)&
$\displaystyle{  \left( \,\overline{R}_{3,1}\right) ^{-2}}$ &
$\displaystyle{  }$ &
$\displaystyle{  \left( \overline{R}_{1,1}\right) ^{-2}}$ &
$\displaystyle{  4 \left( \overline{B}_1 - \underline{B}_1\right) ^{-2}}$ \\ 
& & & &   \\  \hline
& & & &   \\ 

5)&
$({2 \over \pi} \overline{r}_{1,1})^2$ &
$({2 \over \pi} \overline{r}_{3,1})^2$ &
$\displaystyle{  }$ &
$\displaystyle{  }$ \\
& & & &   \\  \hline
& & & &   \\ 

6)&
$({1 \over \pi} \overline{m}_{1})^2$ &
$\displaystyle{  }$ &
$\displaystyle{  }$ &
$\displaystyle{  }$ \\
& & & &   \\  \hline
& & & &   \\ 

7)&
$\displaystyle{  \left({1 \over \pi} \int_0^1 {du \over m_u } \right)^2 }$ &
$\displaystyle{  \left({1 \over \pi} \int_0^1 {du \over r_{3,u} } \right)^2 }$ &
$\displaystyle{  \left({1 \over 2} \int_0^1 {du \over R_{2,u} } \right)^{-2} }$ &
$\displaystyle{  \left({1 \over 2} \int_0^1 {du \over R_{3,u} } \right)^{-2} }$ \\ 
& & & &   \\  \hline
& & & &   \\ 
8)& 
$\displaystyle{  {\tau_1 \over 4 ( \overline{R}_{1,\tau_1} )^2 }  }$ & 
$\displaystyle{  {\tau_1 \over ( \overline{B}_{\tau_1} - \underline{B}_{\tau_1})^2 }  }$ & 
$\displaystyle{   {4 \over \tau_1^2 } \,\int_0^{\tau_1} B_t^2 dt  }$ & 
 \\
& & & &   \\  
\hline
\end{tabular}
\label{table11}
\end{table}

\subsection{Squared Bessel processes (Row 2)}
\label{sqbes}

For $d = 1,2, \ldots$ the squared Bessel process $R_d^2$ is by definition
the sum of $d$ independent copies of $R_1^2 = B^2$, the square of a
one-dimensional Brownian motion $B$, and a similar remark applies to
the squared Bessel bridge $r_d^2$. 
Following \Lev\  \cite{levy51area,lev48}, let us expand the Brownian motion
$(\bmB_t, 0 \le t \le 1)$ or the Brownian bridge $(b_t, 0 \le t \le 1)$ in a 
Fourier series. For example, the standard Brownian bridge $b$
can be represented as
$$
b_u =
\sum_{n=1}^{\infty} {\sqrt 2\over \pi} {Z_n \over n } \sin (\pi n u) ~~~~~~(0 \le u \le 1)$$
 where the $Z_n$ for $n=1,2,\ldots$ are independent standard 
normal random variables, so $E(Z_n) = 0$ and $E(Z_n^2) = 1$ for all $n$.
Parseval's theorem then gives
$$\int_0^1 b_u^2 du =
\sum_{n=0}^{\infty}{Z_n^2 \over\pi^2 n^2}
$$
so the random variable 
$\int_0^1 b_u^2 du$ appears as a quadratic form in 
the  normal variables $Z_n$.
It is elementary and well known that $Z_n^2 \ed 2 \gamma_{1/2}$
for $\gamma_{1/2}$ with gamma$(\hf)$ distribution as in \re{gamhdef}
and \re{gamhlt} for $h = \hf$.
Thus
$$E[\exp (-\lambda Z_n^2)]=(1+2\lambda)^{-1/2}$$ 
and
$$
E\left[\exp \left(-\lambda\int_0^1 b_u^2 du\right)\right]=\prod _{n=1}^{\infty}
\left(1+{2\lambda\over \pi^2 n^2}\right)^{-1/2}
=\left({\sqrt{2\lambda}\over \sinh\sqrt{2\lambda}} \right)^{1/2}
$$
by another application of  Euler's formula  \re{esinh}.
Taking two and four independent copies respectively gives the first two entries
of Row 2. The other entries of this row are obtained by similar
considerations for unconditioned Bessel processes.

Watson \cite{watson61} found that
\eq
\lb{wat}
\int_0^1 \left(b_t - \int_0^1 b_u du \right)^2 dt \ed \qr \TT_{1} .
\en
Shi-Yor \cite{shiyor97} give a proof of \re{wat} with the help of a
space-time transformation of the Brownian bridge. See also
\cite[p. 18-19]{yor92z}, \cite[p. 126-127]{yor97z} and papers cited
there for more general results in this vein.
In particular, we mention a variant of \re{wat} for $B$ instead of $b$, 
which can be obtained as a consequence of a stochastic Fubini theorem 
\cite[p. 21-22]{yor92z}:
\eq
\lb{wat1}
\int_0^1 \left(\bmB_t - \int_0^1 \bmB_u du \right)^2 dt \ed \int_0^1 b_u^2 du \ed \Sigma_{1/2} .
\en
As remarked by Watson, it is a very surprising consequence of \re{wat}
and 
\re{rw1} that 
\eq
\lb{wato}
\int_0^1 \left(b_t - \int_0^1 b_u du \right)^2 dt \ed \pi^{-2}
\max_{0 \le t \le 1} b_t^2 .
\en
As pointed out by Chung \cite{chu76}, the identities in distribution
\re{gngi1} and \re{exclim}, where $\TT_2$ is the sum
of two independent copies of $\TT_1$, imply that the distribution of
$(\max_{0 \le t \le 1} b_t - \min_{0 \le t \le 1} b_t )^2$ is that of
the sum of two independent copies of $\max_{0 \le t \le 1} b_t^2$.
In a similar vein, the first column of Table \tabt\ shows that
the distribution of $4 \pi^{-2} \max_{0 \le t \le 1} b_t^2$ is in turn that of the sum of 
two independent copies of $\int_0^1 b_t ^2 dt \ed \TT_{1/2}$.
There is still no explanation of these coincidences in terms of 
any kind of transformation or decomposition of Brownian paths, or any
combinatorial argument involving lattice paths,
though such methods have proved effective in explaining and generalizing
numerous other coincidences involving the distributions of
$\sigh$ and $\sighp$ for various $h >0$. Vervaat's explanation \re{verv}
of the identity in law between the range of the bridge and the maximum
of the excursion provides one example of this.
Similarly, \re{verv} and \re{wat} imply that
\eq
\lb{watver}
\int_0^1 \left(e_t - \int_0^1 e_u du \right)^2 dt \ed   \qr \sigo .
\en

\subsection{First passage times (Row 3)}
\label{fpt}
It is known \cite{ct62,im65,kent78} that by solving an appropriate
Sturm-Liouville equation, for $\la >0$
$$
E \exp ( - \la T_1(R_d) ) = { (\sqrt{ 2 \la } )^\nu \over 2^\nu \Gamma(\nu +1) I_\nu(\sqrt{2 \la } )}
= \prod_{n = 1}^\infty \left( 1 + {2 \la \over j_{\nu,n} ^2} \right)^{-1}
$$
where $\nu:= (d-2)/2$ with $I_\nu$ the usual modified 
Bessel function, related to 
$J_\nu$ by $(ix)^\nu/I_\nu(ix) = x^{\nu} / J_\nu(x)$,
and 
$j_{\nu,1} < j_{\nu,2} < \cdots$ is the increasing sequence of positive zeros
of $J_\nu$. 
That is to say,
\eq
\lb{specrep}
T_1(R_d) \ed \sum_{n = 1}^\infty {2 \eps_n \over j_{\nu,n}^2 }
\en
where the $\eps_n$ are independent standard exponential variables.
See also Kent \cite{kent80e,kent82s}, and literature cited there,
for more about this {\em spectral decomposition} of $T_1(X)$,
which can be formulated for a much more general one-dimensional
diffusion $X$ instead of $X = R_d$.
The results of Row 3, that
$$
T_1(R_1) \ed \sigop \mbox{ and } T_1(R_3) \ed \sigo 
$$
are the particular cases $d = 1$ and
$d =3$ of \re{specrep}, corresponding to $\nu = \pm 1/2$,
when $I_\nu$ and $J_\nu$ can be expressed in terms of hyperbolic and
trigonometric functions. 
In particular, $j_{-1/2, n} = (n-\hf ) \pi$ and
$j_{1/2, n} = n \pi$ are the $n$th positive zeros of the cosine and 
sine functions respectively.
Comparison of Rows 2 and 3 reveals the identities
$$
T_1(R_1) \ed \int_0^1 R_{2,u}^2 du 
~~~\mbox{and}~~~
T_1(R_3) \ed \int_0^1 r_{2,u}^2 du  .
$$
As pointed out by Williams \cite{wil70,wil74,wil79},
these
remarkable coincidences in distribution are the simplest case $g(u) = 1$ of the identities in law
\eq
\lb{rk1}
\int_0^{T_1(R_1)} g( 1 - R_{1,t})  dt  \ed \int_0^1 R_{2,u}^2 g(u) du 
\en
and
\eq
\lb{rk11}
\int_0^{T_1(R_3)} g(R_{3,t})  dt  \ed \int_0^1 r_{2,u}^2 g(u) du 
\en
where the two Laplace transforms involved are again determined by the
solutions of a Sturm-Liouville equation \cite{py82}, \cite[Ch. XI]{ry99}.
Let $L_t^x(R_d), t \ge 0, x \in \reals)$ and
$L_t^x(r_d), 0 \le t \le 1, x \in \reals)$ be the
local time processes of $R_d$ and $r_d$ defined by the occupation
density formula \re{ltdef} with $B$ replaced by $R_d$ or $r_d$.
Granted existence of local time processes for $R_d$ and $r_d$,
the identities \re{rk1} and \re{rk11} are an expression of the Ray-Knight 
theorems \cite[Ch. XI, \S 2]{ry99} that
\eq
\lb{rk2}
(L_{T_1(R_1)}^{1-u}(R_1), 0 \le u \le 1) \ed (R_{2,u}^2, 0 \le u \le 1)
\en
and
\eq
\lb{rk3}
(L_{T_1(R_3)}^{u}(R_3), 0 \le u \le 1) \ed (r_{2,u}^2, 0 \le u \le 1) .
\en
The next section gives 
an interpretation of the variable
$T_1 (R_3) + T_1(\hat{R}_3)$ appearing in column 2 in terms
of Brownian excursions.

\subsection{Maxima and the agreement formula (Rows 4 and 5)}
\label{max}

The entries in Row 4 are equivalent to corresponding entries in Row 3 by
application to $X = R_d$ and $X = \overline{B} - \underline{B}$
of the elementary identity 
\eq
\lb{inverse}
(\overline{X}_{1} )^{-2} \ed T_1(X) 
\en
which is valid for any process $X$ with continuous paths which
satisfies the Brownian scaling identity \re{bscale}, because 

$$
P((\overline{X}_{1}) ^{-2} >t) =
P ( \overline{X}_{1} < t^{-\hf}) = P ( \overline{X}_{t} < 1 )  =
P(T_1 (X)> t ) .
$$
The first entry of Row 5, with
$\overline{r}_{1,1}:= \max_{0 \le u \le 1} |b_u|$,
is read from \re{rw1}.
The second entry of Row 5, involving the maximum $\overline{r}_{3,1}$ of a 
three-dimensional Bessel bridge $(r_{3,u}, 0 \le u \le 1)$, is read
from the work of Gikhman \cite{gik} and Kiefer \cite{kie58}, 
who found a formula for $P(\overline{r}_{d,1} \le x)$ for arbitrary 
$d = 1,2, \ldots$. See also \cite{py97max}. 
This result involving
$\overline{r}_{3,1}$ may be regarded as a consequence of the previous
identification \re{exclim} of the law of $\overline{e}_{1}$, the 
maximum of a standard
Brownian excursion, and the identity in law 
$ \overline{e}_{1} \ed \overline{r}_{3,1} $
implied by the remarkable
result of \Lev-Williams \cite{lev48,wil70}, that
\eq
\lb{levwil}
(e_t, 0 \le t \le 1) \ed (r_{3,t}, 0 \le t \le 1) .
\en
Another consequence of the scaling properties of Bessel processes is 
  provided by the following absolute
continuity relation between the law of 
$(\overline{r}_{d,1})^{-2}$ and the law of 
$$
\sigtd:= T_1 (R_d) + T_1(\hat{R}_d)
$$ 
for general $d>0$.
This result, obtained in \cite{by87,py93c,bi90wuhan,py95agr}, we call the
{\em agreement formula}: for every non-negative Borel function $g$
\eq
\lb{agr1}
E \left[ g(
(\overline{r}_{d,1})^{-2} \right)] = 
C_d E \left[ \sigtd ^{\nu}  \, g(\sigtd) \right]
\en
where 
$C_d := 2^{(d-2)/2}\Gamma(d/2)$.
In \cite{py95agr} the agreement formula was presented 
as the specialization to Bessel processes of a 
general result for one-dimensional diffusions. 
As explained in \cite{by87,wzeta,py95agr}, the 
agreement formula follows from the fact that a certain
$\sigma$-finite measure on the space of continuous non-negative paths with 
finite lifetimes can be explicitly disintegrated in two different ways, 
according to the lifetime, or according to the value of the maximum.

Note from \re{specrep} that $\Sigma_{2,3} \ed \sigt$ and 
$\Sigma_{2,1} \ed \sigtp$.
For $d=3$ formula \re{agr1} gives for all non-negative Borel 
functions $g$
\eq
\lb{agree}
E \left[ g( \overline{r}_{3,1} ) \right]  = \rtop E \left[ \sqrt{ \sigt } ~ g( 1/\sqrt{\sigt}) \right] .
\en
In view of \re{agree}, the symmetry property 
\re{gheq} of the common distribution of $Y$ and 
$\sqrt{ { \pi \over 2 } \sigt }$, which expresses the functional equations for 
$\xi$ and $\theta$, can be recast as the following identity of
Chung \cite{chu76}, which appears in the second column of Table \tabt:
\eq
\lb{chungid}
\left( \mbox{${2 \over \pi }$} \,\overline{r}_{3,1} \right) ^2 \ed \sigt .
\en
As another application of \re{agr1}, we note that for $d = 1$ this formula
shows that the reciprocal relation between the laws of $\sigo$ and
$\sigtp$ discussed in Section \ref{infdiv} is equivalent to 
the equality in distribution of \re{rw1}, that is
\eq
\lb{chungid1}
\left( \mbox{${2 \over \pi }$} \,\overline{r}_{1,1} \right) ^2 \ed \sigo .
\en
We do not know of any path transformation leading to a non-computational
proof of \re{chungid} or \re{chungid1}.

\subsection{Further entries.}
The distributions of $T_1(R_d)$ and 
$T_1(R_d) + T_1 (\hat{R}_d)$ for $d =1,3$ shared by the columns of Table \tabt, 
also arise naturally from a number of other constructions involving
Brownian motion and Bessel processes.
Alili \cite{alili96ib} found the remarkable result that 
\eq
\lb{alili}
{ \mu^2 \over \pi^2 } \left[ \left(\int_0^1 \coth ( \mu r_{3,u}) du \right)^2 -1 \right] \ed \sigt
\mbox{ for all } 
\mu  \ne 0.
\en
As a check, the almost sure limit of the left side of \re{alili}
as $\mu \te 0$ is the variable 
$\pi^{-2} (\int_0^1 r_{3,u}^{-1} du )^2$ in the second column 
of Row 7. As shown by Alili-Donati-Yor \cite{admy97}, consideration
of \re{alili} as $\mu \te \infty$ shows that
%the left side of \re{alili} converges in law to the sum of two independent
%copies of
\eq
\lb{admy}
{4 \over \pi^2} \int_0^\infty{ dt \over \exp ( R_{3,t} ) - 1}  \ed \sigo
\en
Other results are the identity of Ciesielski-Taylor \cite{ct62}
according to which
\eq
\lb{ct} 
\int_0^\infty 1( R_{d+2,t} \le 1 ) dt \ed  T_1(R_d)
\en
which for $d = 1$ and $d =3$ provides further entries for the table.
See also \cite{kt78}, \cite{jp95cyc}, \cite[p.97-98, Ch. 7]{yor92z},
\cite[p. 132-133]{yor97z}, \cite{py93lt} for still more 
functionals of 
Brownian motion whose Laplace transforms can be expressed in terms of 
hyperbolic functions.
\section{Renormalization of the series $\sum n^{-s}$.}
\label{renorm}
\subsection{Statement of the result}
\label{Theorem1}
The expansion of $\sigo$ as an infinite series \re{sighdef} 
suggests that we use
partial sums in order to approximate its Mellin transform. As we shall see, this
yields an interesting approximation of the Riemann zeta function.
Consider again the relationship \re{xifn} between $\zeta$ and $\xi$,
which allows the definition of $\zeta(s)$ for $s \ne 1$ despite the
lack of convergence of \re{zeta1} for $\Re s \le 1$. There are a number of known
ways to remedy this lack of convergence, some of which are discussed in Section
\ref{summation}.  
One possibility is to look for an array 
of coefficients $(a_{n,N}, 1 \le n \le N)$ such that the functions
\eq \lb{chi}
\kappa_N(s):=\sum_{n=1}^N{a_{n,N}\over n^s}
\en
converge as $N \te \infty$, for all values of $s$.
For fixed $N$ there are $N$ degrees of freedom in the choice of the 
coefficients, so we can enforce the conditions $\kappa_N(s)=\zeta(s)$ at 
$N$ choices of $s$, and it is natural to choose the points $0$, where 
$\zeta(0)=-{1\over 2}$, and $-2,-4,-6,\ldots, -2(N-1)$ where $\zeta$
vanishes. 
It is easily checked that this makes
\eq
\lb{anN}
a_{n,N}= {(-N)(1-N)(2-N)\ldots (n-1-N)\over
 (N+1)(N+2)\ldots (N+n)}=(-1)^n{\pmatrix{2N\cr
 N-n}\over \pmatrix{2N \cr N}} .
\en
Note that for each fixed $n$
$$
a_{n,N} \te (-1)^n \mbox{  as } N\te \infty
$$
and recall that 
$$
(2^{1-s}-1)\zeta(s)=\sum_{n=1}^\infty{(-1)^n\over n^s} \mbox{   for } \Re s>1
$$ 
extends to an entire function of $s$. One has

\begin{thm}{Thm1}
For $\eps_i, 1 \le i \le N$ independent standard exponential variables,
and $\Re s > - 2N$ 
\eq
\lb{chiN}
E \left[ \left( {2 \over \pi^2} \sum_{n = 1} ^N { { \eps_n }\over n^2}
\right) ^{s/2} \right] = - s {2^{s/2} \over \pi^s }
\Gamma(\hfs) \, \sum_{n=1}^N{a_{n,N}\over n^s}
\en
where the $a_{n,N}$ are defined by \rem{anN}, and
\eq
\lb{zappr}
\sum_{n=1}^N{a_{n,N}\over n^s} \te (2^{1-s}-1)\zeta(s) \mbox{ as  }N\te\infty
\en
uniformly on every compact subset of $\complex$.
\end{thm}

The proof of Theorem \ref{Thm1} occupies the next two sections.
\subsection{On sums of independent exponential random variables}
\label{sumex}

Let $(\expe_n;n\geq 1)$ be a sequence of independent identically distributed
random variables, with the standard exponential distribution 
$P(\expe_n \ge x) = e^{-x}$,$x \ge 0$.
Let $(a_n;n\geq 1)$ be a sequence of positive real numbers, such that 
$
\sum_{n=1}^{\infty}
 a_n<\infty
 $,
  then the series 
 $
 X=\sum_{n=1}^{\infty}
  a_n \expe_n
  $
   converges almost surely, and
in every $L^p$ space, for $1\leq p<\infty$.
\begin{lmm}{conv}
Let  $X=\sum_{n=1}^{\infty} a_n \expe_n$ as above, and let 
$X_N= \sum_{n=1}^N a_n \expe_n$ be the partial sums, then
for every real $x$ one has $E[X^x]<\infty$, and
$E[X_{N}^x]<\infty $ for $x>-N$. Furthermore one has 
$$E[X_N^s]\te
E[X^s]~~~~~ \mbox{as }{N\te\infty} 
$$ uniformly with respect to $s$ on each compact subset of $\complex$.
\end{lmm}
\proof
We have already seen that $E[X^x]<\infty$ if $x\geq 0$. Let us prove that
$E[X_{N}^x]<\infty $ for $0>x>-N$. Let $b_N=\min(a_n;n\leq N)$, then
$X_N\geq b_NY_N=b_N\sum_{n=1}^N\expe_n$. But  $Y_N$ has a gamma distribution,
 with density ${1\over\Gamma(N)}t^{N-1}e^{-t}$ at $t >0$, so that $E[Y_N^x]<\infty$ for
  $0>x>-N$ and thus $E[X_N^x]<\infty$.
 The assertion for $X$ follows from $X\geq X_N$. It remains to check the
 uniform convergence. If $\Re s\in [-A,+A]$, then
 \begin{eqnarray*}
\lb{aaaa}
\vert X_N^s-X^s\vert&= & \left| \int_{X_N}^Xsy^{s-1}dy\right|\\
 & \leq &\vert s\vert\, (X-X_N)(X_N^{-A-1} \vee X^{A-1})
\end{eqnarray*}
and the required uniform convergence as $N \te \infty$ is now evident by
application of the Cauchy-Schwarz inequality.
\endpf

We now compute the
Mellin transform of the distribution of $\sum_{n=1}^N a_n\expe_n$, assuming
that the $a_n$ are all distinct and strictly positive.
\begin{lmm}{mellinexp}
With the above notations, and 
$\Pin:= \prod_{j\ne n, 1 \le j \le N}\left(1-{a_j\over a_n}\right)$, for $\Re s>-N$
\eq
\lb{emel}
E\left[ \,\left( \sum_{n = 1}^N a_n \expe_n \right)^s \,\, \right]  
= \Gamma(s+1) \,
\sum_{n=1}^N {a_n^{s}\over \Pin} 
\en
where the right side of \re{emel} is defined by continuity for
$s=-1,-2,\ldots, -N+1$.
\end{lmm}
\proof
The partial fraction expansion of the Laplace transform
\eq
\lb{ltpart}
E\left[ \exp \left( - \la \sum_{n = 1}^N a_n \expe_n \right) \,\, \right]  
= 
\prod_{n=1}^N {1 \over (1 + \la a_n)} =
\sum_{n=1}^N {1 \over \Pin} {1 \over (1 + \la a_n)}
\en
implies that for every non-negative measurable function $g$ such that
$E[ g(a_n \eps_1)]$ is finite for every $n$
\eq
\lb{giden}
E\left[ g \left( \sum_{n = 1}^N a_n \expe_n \right) \,\, \right]  
= \sum_{n=1}^N {1 \over \Pin} E [ g( a_n \eps_1) ] .
\en
For $g(x) = x^s$ this gives \re{emel}, first for real $s > -1$, then
also for $\Re s > -N$ since the previous lemma shows that the left 
side is analytic in this domain, and the right side is evidently 
meromorphic in this domain.
\endpf

Note the implication of the above argument that the sum on the right side
of \re{emel} must vanish at $s=-1,-2,\ldots, -N+1$.
\subsection{Proof of Theorem 3}
\label{prooft1}

This is obtained by applying the preceding results in the particular
case $a_n = n^{-2}$.
By application of Lemmas \ref{conv} and \ref{mellinexp}, 
and the formula for $E(\sigo^s)$ in Table 1, found in 
\cite[(86)]{py97max}, the conclusions of Theorem \ref{Thm1} hold for 
\eq
\lb{anNp}
- 2 a_{n,N} =  {1 \over \Pin} = 
{ \prod_{j \ne n} (j^2 ) \over \prod_{j \ne n} (j^2  -n^2) }
\en
where both products are over $j$ with  $1 \le j \le N$ and
$j \ne n$.
The product in the numerator is $(N!/n)^2$ while writing
$j^2  -n^2 = (j-n)(j+n)$ allows the product in the denominator to 
be simplified to $(-1)^{n-1}(N+n)!(N-n)!/(2 n^2)$. Thus the expression
for $a_{n,N}$ can be simplified to \re{anNp}.
\endpf

\subsection{The case of the $L_{\chi_4}$ function}
The following result can be obtained similarly, with the help of the formula 
%\eq
%\lb{y11b}
%E \left[ \left(  \sigop \right) ^{s/2} \right] = 
%\Gamma(s+1) 2^{s+1} \left( 2 \over \pi \right)^{2s+1} L_{\chi_4}(2 s + 1)
%\en
for $E ((\sigop )^s)$ in
terms of $L_{\chi_4}$ defined by \re{oddser},
given in Table 1 of Section \ref{infdiv}. 
\begin{thm}{zeta}
Let 
$$
 L_{\chi_4}^{(N)}(s):=\sum_{n=0}^{N-1}(-1)^{n}{\pmatrix{2N-1\cr N-n-1}\over
 \pmatrix{2N-1\cr N-1}}(2n+1)^{-s}.
$$
Then 
${\pmatrix {2N-1\cr N-n-1}\over\pmatrix{2N-1\cr N}} \to 1$
as $N\to \infty$;  for each
$N$,
one has
$L_{\chi_4}^{(N)}(1-2k)=0$ for $k=1,2,\ldots N-1$, and $L_{\chi_4}^{(N)}
(s)\to L_{\chi_4}(s)$
uniformly on every compact of $\complex$.
\end{thm}
We note that it is also possible to use formulae \re{y1} and \re{sigy} 
to provide another approximation of $\zeta$. We leave the computation to the
interested reader.
\subsection{Comparison with other summation methods}
\label{summation}
Perhaps 
the simplest way to renormalize the series \re{zeta1} is given by the 
classical formula 
\eq
\zeta(s)=\lim_{N\to\infty}\left(\sum_{n=1}^Nn^{-s}-{N^{1-s}-1\over
1-s}\right)-{1\over 1-s}~~~~(\Re(s)>0)
\en
Related methods are provided by the 
approximate functional equation and Riemann-Siegel formula, which
 are powerful tools in
deriving results on the behaviour of the zeta function in the critical strip.
See e.g. Ch. 7 of Edwards \cite{edwards74} for a detailed discussion.

It is also known \cite{landau07,hardy49div}
that the  series $\sum_1^{\infty}{(-1)^n\over n^s}$
 is Abel summable for all values of $s\in\complex$, meaning that as $z\to 1$ in the
unit disk,  
$$
\sum_1^{\infty}{(-z)^n\over n^s}\to (2^{1-s}-1)\zeta(s).
$$
The Lerch zeta function
$$
\Phi(x,a,s)=\sum _{n=0}^{\infty}{e^{2i\pi n x}\over (n+a)^s}~~~~~(x\in\real,\,
0<a\leq 1,\, \Re(s)>1)
$$
is known to have analytic continuation to $s\in\complex$, with a pole at $s=1$
 for $a=1$,
$x\in\ints$. 
This allows us to
sketch another proof of Theorem \ref{Thm1}. The formula
$$
\kappa_N(s)={2^{2N}}{\pmatrix{2N\cr N}}^{-1} \,  \int_0^1 (\sin(\pi
x))^{2N}\Phi(x,1,s)\,dx
$$
is easily checked using \re{chi}-\re{anN} for $\Re(s)>1$, 
and extended by analytic continuation to all values of
$s\in\complex$. Convergence of $\kappa_N(s)$ towards
$(2^{1-s}-1)\zeta(s)$ then follows from continuity properties of the 
Lerch zeta function in the variable $x$, and the fact that
${2^{2N}} {\pmatrix{2N\cr N}} ^{-1} \, (\sin(\pi x))^{2N} \,dx\to \delta_{1/2}$
 weakly as $N\to\infty$.

Finally, we note that
J. Sondow \cite{sondow94} has shown that Euler's summation method yields
the following series, uniformly convergent on every compact of $\complex$:
\eq
\lb{sondow}
(1-2^{1-s})\zeta(s)=
\sum_{j=0}^{\infty}{1-\pmatrix{j\cr 1}2^{-s}+
\ldots+(-1)^j\pmatrix{j\cr j}(j+1)^{-s}
\over 2^{j+1}} .
\en
Furthermore the  sum of the first $N$ terms of this  series gives
 the exact values of $\zeta$ at $0,-1,-2,\ldots, -N+1$, so we can rewrite
 the partial sum in \re{sondow} as
 $$
 \rho_N(s)=\sum_1^{N}{b_{n,N}\over n^s}
 $$
 where the $b_{n,N}$ are completely determined by 
 $\rho_N(-j)=\zeta(-j)$ for $j=0,1,\ldots, N-1$. 
Compare with the discussion between \re{chi} and \re{anN} to see 
the close parallel between \re{zappr} and \re{sondow}.

\section{Final remarks}
\label{sfc}
\subsection{Hurwitz's zeta function and Dirichlet's $L$-functions}
\label{hur}
Row 3 of Table \tabt\ involves hitting times of Bessel processes of dimension 1
and 3, started from 0. If the Bessel process does not start from zero, we still
have an interesting formula for the Mellin transform of the hitting time,
expressed now in terms of the Hurwitz zeta function. Specifically, one has
\eq
\lb{hita}
E[e^{-\la T_1^a(R_3)}]={\sinh(a\sqrt{2\la})\over a\sinh(\sqrt{2\la})}
 ~~~~~E[e^{-\la T_1^a(R_1)}]={\cosh(a\sqrt{2\la})\over\cosh(\sqrt{2\la})}
 \en
 where $T_1^a$ denotes the hitting time of 1, starting from $a\in]0,1[$,
  of the
 corresponding Bessel process. Expanding the denominator we get 
 $$
 {\sinh(a\sqrt{2\la})\over a\sinh(\sqrt{2\la})}={1\over a}\sum_{n=0}^{\infty}
 e^{-(2n+1-a)\sqrt{2\la}}-e^{-(2n+1+a)\sqrt{2\la}}
 $$
 Inverting the Laplace transform yields the density of the distribution of 
 $T_1^a(R_3)$
 $$
 {1\over a \sqrt{2\pi t^3}}\sum_{n=0}^{\infty}
(2n+1-a) e^{-(2n+1-a)^2/(2t)}-(2n+1+a)e^{-(2n+1+a)^2/(2t)}
 $$
 Taking the Mellin transform we get
 $$
 E[(T_1^a(R_3))^{s/2}]={\Gamma({s-1\over 2})\over a 2^{s/2}}
\left(\zeta(s,{1-a\over 2})-
 \zeta(s,{1+a\over 2})\right)~~~~~(\Re(s)>1)
 $$
where $\zeta(s,x)=\sum_{n=0}^{\infty}(n+x)^{-s}$ is the Hurwitz zeta function.
This identity extends by analytic continuation to all $s\in\complex$.
A similar expression exists for $T_1^a(R_1)$.

One can use the product expansion for $\sinh$ in order to give an
approximation of 
$\zeta(s,u)- \zeta(s,1-u); u\in]0,1[$. For it is easy to see
that $\prod_{n=1}^N\left(1+{2a^2\la\over n^2}\right)\left(1+{2\la\over n^2}\right)^{-1}$ is the Laplace
transform of a probability distribution on $[0,\infty[$,
and that  this probability distribution converges towards that
of $T_1^a(R_3)$, in such a way that there is a result similar to 
Theorem \ref{Thm1}.

The 
Hurwitz zeta function can be used to construct Dirichlet $L$-functions by linear
combinations. However, direct probabilistic interpretations of general Dirichlet
$L$-functions, in the spirit of what we did in Section \ref{brown} do not seem
to exist. More precisely, 
let $\chi$ be a primitive character modulo $N$, and let
\eq
\lb{thetachie}
\theta_{\chi}(t)=\sum_{n=-\infty}^{+\infty}n ^{\epsilon} 
\chi(n)e^{-\pi n^2t}
\en
%\eq
%\lb{thetachio}
%\theta_{\chi}(t)=\sum_{n=-\infty}^{+\infty}n\chi(n)e^{-\pi n^2t}~~~~
%\hbox{if $\chi$ is odd, i.e. $\chi(-1)=-1$.}
%\en
where $\epsilon=0$ or $1$ according to whether
$\chi$ is even or odd, so $\chi(-1)=(-1)^{\epsilon}$.
These functions satisfy the functional equation 
\eq
\lb{fethetachi}
\theta_{\chi}(t)={(-i)^{\epsilon}\tau(\chi)\over
N^{1+\epsilon}t^{\epsilon+1/2}}\theta_{\bar\chi}\left({1\over N^2t}\right)
\en
where $\tau(\chi)$ is a Gauss sum.
Taking a Mellin transform, 
this yields the analytic continuation
and functional equation for the associated Dirichlet $L$-function
$$
L_\chi(s):=\sum_{n=1}^{\infty}{\chi(n)\over n^s},
$$
namely
\eq
\lb{Lfe}
\Lambda(s,\chi)=(-i)^{\epsilon}\tau(\chi)N^{-s}\Lambda(1-s,\bar\chi)
\en
where 
$$
\Lambda(s,\chi)=\pi^{{-(s+\epsilon) \over 2}}\Gamma\left({s+\epsilon\over 2}\right)\,
L_\chi (s) .
$$
See \cite{dav52} or \cite[\S 1.1]{bump97} for the classical derivations
of these results.
 For general real $\chi$, there does not seem to be any simple
probabilistic interpretation of $\theta_{\chi}(t)$.
In particular, this function is not necessarily positive for all $t >0$.
This can be seen as follows. We choose an odd character $\chi$ 
 (the case of even characters is similar), and 
 compute the Laplace transform of $\theta_{\chi}$ using
 \re{fethetachi}
 $$
 \int_0^{\infty}e^{-\la t}\theta_{\chi}(t)\,dt=\sum_{n=-\infty}^{+\infty}
 \int_0^{\infty}{n\chi(n)\over N^{3/2}t^{3/2}}e^{-\la
 t}e^{-\pi n^2/(N^2t)}\,dt=N^{-1/2}\sum_{n=1}^{\infty}
 \chi(n)e^{-{2n\over N}\sqrt{\pi\la}}
 $$
 Using the periodicity of $\chi$,
  this equals
 \eq
 \lb{lapchi}
 {N^{-1/2}
 \sum_{n=1}^{N-1}\chi(n)e^{-{2n\over N}\sqrt{\pi\la}}\over 1-e^{-2\sqrt{\pi\la}}}
 ={N^{-1/2}\sum_{n=1}^{(N-1)/2}\chi(n)\sinh({N-2n\over N}\sqrt{\pi\la})\over
 \sinh(\sqrt{\pi\la})}
 \en
 For small values of $N$, and $\chi$ a real odd character modulo $N$,
  one can see by inspection that this indeed is the
 Laplace transform of a positive function, hence by uniqueness of Laplace
 transform, $\theta_{\chi}$ is  positive on the real line.
 However \Pol\ \cite{polya19} exhibited 
   an infinite number of primes $p$ such that for the
 quadratic character modulo $p$ the polynomial $Q_p(x):= \sum_{n=1}^{p-1} x^n\chi(n)$ 
 takes negative
 values somewhere  on $[0,1]$. 
% These primes  satisfy
% $\chi(2)=\chi(3)=\chi(5)=\chi(7)=\chi(11)=\chi(13)=-1$, and their existence is 
% secured by invoking  Gauss quadratic
% reciprocity, and Dirichlet's theorem on primes in arithmetic progressions.
In particular, for $p = 43$ we find $Q_{43}(3/4) \approx -0.0075$.
For such quadratic Dirichlet characters, 
  the Laplace transform above also 
 takes negative
 values, which implies that $\theta_{\chi}$ does not stay positive on $[0,\infty[$.
 We note that if $\theta_{\chi}>0$ on $]0,\infty[$ then obviously
  its Mellin transform has no zero on the real line, and hence 
  the corresponding $L$-function has no Siegel zeros.

\subsection{Other probabilistic aspects of Riemann's zeta function}
\label{ref}
It is outside the scope of this paper, and beyond
 the competence of its authors, 
to discuss at length the theory of the Riemann zeta function. 
But we mention in this final Section some other 
works relating the zeta function to probability theory. 

The work of \Polya\
 has played a significant role in the proof of the
Lee-Yang theorem in statistical mechanics:
 see the discussion in \cite[pages 424-426]{polya74}.
Other connections between Riemann zeta function 
and statistical mechanics appear in 
 Bost and Connes \cite{bosco95} and in Knauf \cite{knauf97}.

The Euler product for the Riemann zeta 
function is interpreted probabilistically
in Golomb \cite{golomb70} and  Nanopoulos \cite{nano4},
via the independence of various prime factors
when choosing a positive integer according to the distribution with
probability at $n$ equal to $\zeta(s)^{-1} n^{-s}$ for some $s >1$.
See also Chung \cite[p. 247]{chung74} and \cite{sd88}.

It is an old idea of Denjoy \cite{denj64} that the partial sums of the 
M\"obius function should behave like a random walk (the law of iterated logarithm would imply Riemann hypothesis). 

Fascinating connections between zeros of the Riemann zeta function (and other
$L$-functions) and eigenvalue distribution of random matrices 
are currently under intense scrutiny. See 
Odlyzko \cite{odlyzko87} and Katz and Sarnak \cite{ksarnak99}.

Finally let us mention a few other recent references where the Riemann zeta
function makes some appearances: Asmussen, Glynn and Pitman \cite{agp93},
Joshi and Chakraborty \cite{jc96}, Chang and Peres \cite{changp97}.
  
%%to recompile the references
%\input{0.bibcall}

%%to input the references

\end{document}